\documentclass[12pt]{article}
\usepackage{amsmath}
\usepackage{amsfonts}
\usepackage{latexsym}
\usepackage{amssymb}
\usepackage{a4wide}
\usepackage[all]{xy}
\newtheorem{thm}{Theorem}[section]
\newtheorem{la}[thm]{Lemma}
\newtheorem{Defn}[thm]{Definition}

\newtheorem{Exam}[thm]{Example}
\newtheorem{Remark}[thm]{Remark}

\newtheorem{prop}[thm]{Proposition}
\newtheorem{cor}[thm]{Corollary}
\newenvironment{defn}{\begin{Defn}\rm}{\end{Defn}}
\newenvironment{ex}{\begin{Exam}\rm}{\end{Exam}}
\newenvironment{rem}{\begin{Remark}\rm}{\end{Remark}}
\newenvironment{proof}{{\noindent\bf Proof.}}%
                  {\nopagebreak\hspace*{\fill}$\Box$\medskip\medskip\par}
\newenvironment{Proof}[1]{{\noindent\bf Proof #1.}}%
                  {\nopagebreak\hspace*{\fill}$\Box$\medskip\medskip\par}
\newcommand{\Punkt}{\nopagebreak\hspace*{\fill}$\Box$}
\newcommand{\wb}{\overline}
\renewcommand{\phi}{\varphi}
\newcommand{\ra}{\rightarrow}
\newcommand{\xra}{\xrightarrow}
\newcommand{\tri}{\triangle}
\newcommand{\mc}[1]{\mathcal{#1}}

\newcommand{\impl}{\Rightarrow}

\newcommand{\mto}{\mapsto}

\newcommand{\ve}{\varepsilon}

\newcommand{\isom}{\cong}

\newcommand{\N}{{\mathbb N}}

\newcommand{\F}{{\mathbb F}}

\newcommand{\C}{{\mathbb C}}

\newcommand{\T}{{\mathbb T}}
\newcommand{\A}{{\mathbb A}}

\newcommand{\SU}{{\rm SU}}

\newcommand{\Sp}{{\rm Sp}}

\newcommand{\I}{{\mathbb I}}

\newcommand{\cO}{{\mathcal O}}
\newcommand{\cC}{{\mathcal C}}

\newcommand{\cS}{{\mathcal S}}
\newcommand\CC{\mathbb{C}}

\newcommand{\cT}{{\mathcal T}}
\newcommand{\g}{{\mathfrak g}}

\newcommand{\dl}{{\displaystyle \lim_{\longrightarrow}}}
\DeclareMathOperator{\End}{End}

\DeclareMathOperator{\Aut}{Aut}

\newcommand{\take}{\backslash}
\newcommand{\Lg}{\mathfrak g}
\newcommand{\Lh}{\mathfrak h}
\newcommand{\Lk}{\mathfrak k}

\newcommand{\sub}{\subseteq}

\DeclareMathOperator{\GL}{GL}
\DeclareMathOperator{\SL}{SL}

\DeclareMathOperator{\pr}{pr}

\DeclareMathOperator{\Mor}{Mor}

\DeclareMathOperator{\id}{id}

\newcommand{\cB}{{\mathcal B}}
\newcommand{\cA}{{\mathcal A}}

\newcommand{\lb}{\left\{ }
\newcommand{\rb}{\right\} }

\DeclareMathOperator{\Repart}{Re}

\DeclareMathOperator{\ad}{ad}
\newcommand{\pl}{{\displaystyle\lim_{\longleftarrow}}}

\newcommand{\G}{\mathbb{G}}
\newcommand{\TG}{\mathbb{T}\mathbb{G}}
\newcommand{\HTG}{\mathbb{H}\mathbb{T}\mathbb{G}}
\newcommand{\LCG}{\mathbb{L}\mathbb{C}\mathbb{G}}
\newcommand{\LIE}{\mathbb{L}\mathbb{I}\mathbb{E}}

\newcommand{\KOG}{\mathbb{K}\mathbb{O}\mathbb{G}}
\DeclareMathOperator{\ob}{ob}

\mathchardef\tnode="020E 

\def\arc{
  \hbox{\kern -0.15em
  \vbox{\hrule width 3em height 0.6ex depth -0.5 ex}
  \kern -0.33em}}

\def\darc{
  \rlap{\lower0.2ex\arc}{\raise0.2ex\arc}}

\def\tarc{
  {\rlap{\arc}{\rlap{\lower0.4ex\arc}{\raise0.4ex\arc}}}}

\def\stroke#1{
  \kern 0.05em
  \rlap\arc{{\textstyle{#1}}\atop\phantom\arc}
  \kern -0.22em}

\def\dstroke#1{
  \kern 0.05em
  \rlap\darc{{\textstyle{#1}}\atop\phantom\darc}
  \kern -0.22em}

\def\tstroke#1{%
  \kern 0.05em
  \rlap\tarc{{\textstyle{#1}}\atop\phantom\tarc}
  \kern -0.22em}

\def\centerscript#1{
  \setbox0=\hbox{$\tnode$}
  \hbox to \wd0{\hss$\scriptstyle{#1}$\hss}}

\def\node{
  \def\super{}
  \def\sub{}
  \futurelet\next\dolabellednode}

  \let\sp=^
  \let\sb=_

  \def\dolabellednode{%
    \ifx\next\sb\let\next\getsub
    \else
      \ifx\next\sp\let\next\getsuper
      \else\let\next\donode
      \fi
    \fi
    \next}

  \def\getsub_#1{\def\sub{#1}\futurelet\next\dolabellednode}
  \def\getsuper^#1{\def\super{#1}\futurelet\next\dolabellednode}

  \def\donode{%
   \rlap{$\mathop{\phantom\tnode}\limits_{\centerscript{\sub}}^{\centerscript{\super}}$}\tnode}

\def\varcdn{
  \kern -0.03em\vbox{\kern -0.5ex
  \hbox to \wd0{\hss\vrule width 0.04em depth 5.8ex\hss}
  \kern -0.3ex  \hbox{$\tnode$}}}

\begin{document}
\begin{center}
{\Large\bf Final Group Topologies, Kac-Moody Groups and\\[2mm]
Pontryagin Duality}\vspace{4.7mm}\\
{\bf Helge Gl\"{o}ckner,
Ralf Gramlich
and Tobias Hartnick}\vspace{1.2mm}
\end{center}
{\bf Abstract.} We study final group topologies and their
relations to compactness properties.\linebreak
In particular, we are
interested in situations where a colimit or direct limit is
locally compact, a $k_\omega$-space, or locally~$k_\omega$. As a
first application, we show that unitary forms of complex Kac-Moody
groups can be described as the colimit of an amalgam of subgroups
(in the category of Hausdorff topological groups,
and the category of $k_\omega$-groups). Our second
application concerns Pontryagin
duality theory for the classes of
almost metrizable
topological abelian groups, resp.,
locally~$k_\omega$
topological abelian groups, which are dual to each
other. In particular, we explore the relations between countable
projective limits of
almost metrizable abelian groups and countable direct
limits of locally $k_\omega$ abelian groups.%
\vspace{1.9mm}\renewcommand{\thefootnote}{\fnsymbol{footnote}}
\footnotetext{\hspace{-2mm}{\em Classification\/}:
22A05 (main); 
18A30, 
20E06, 
22C05, 
22D35, 
22E15, 
43A40, 
51B25, 
51E24, 
54D50\\ 
\hspace*{4.3mm}{\em Keywords\/}: Final group topology,
colimit, direct limit,
inductive limit, projective limit, inverse limit, amalgam,
locally compact group, $k_\omega$-space,
locally $k_\omega$ space,
$k_\omega$-group,
reflexive group,
almost metrizable group,
\v{C}ech complete group,
Pontryagin duality,
Kac-Moody group, building, twin building, Phan system}
\renewcommand{\thefootnote}{\arabic{footnote}}
\begin{center}
{\bf\Large Introduction}\vspace{-.2mm}
\end{center}
Given a group~$G$ and a family
$(f_i)_{i\in I}$ of maps
$f_i\colon X_i\to G$
from certain topological spaces to~$G$,
there exists a finest group topology
on~$G$ making all of the maps $f_i$
continuous, the so-called
\emph{final group topology}
with respect to the family $(f_i)_{i\in I}$.
Such topologies arise naturally
in connection with colimits
of topological groups
(notably, direct limits),
which carry the final group topology
with respect to the family of limit maps.
Although a final group topology~$\cO$
always exists, it may be quite elusive
in the sense that it may not be clear at all
how one could check whether a given
subset $U\sub G$ belongs to~$\cO$.\\[2.5mm]
For example, consider an ascending sequence
$G_1\sub G_2\sub\cdots$ of topological groups
such that all inclusion maps are continuous homomorphisms.
Then $G:=\bigcup_{n\in \N}G_n$ can be made
the direct limit group $\dl\,G_n$,\vspace{-.7mm}
and one can consider a much more concrete
topology~$\cT$ on~$G$, the \emph{direct limit topology},
defined by declaring $U\sub G$ open if and only if
$U\cap G_n$ is open in~$G_n$ for each $n\in \N$.
But unfortunately,
$\cT$ need not make the group
multiplication $G\times G\to G$ continuous,
in which case~$\cT$
is properly finer than the final group
topology~$\cO$ (see \cite{TSH}).
If this pathology occurs, then the natural continuous bijection
\begin{equation}\label{ctsncts}
\dl\, (G_n\times G_n)
\; \to\;
\dl\, G_n\, \times \, \dl\,G_n
\end{equation}
is not a homeomorphism
(cf.\ \cite{Hir}), and thus the pathology is related to
the non-compatibility of direct products
and direct limit topologies.\\[2.5mm]
However, both pathologies disappear
if each $G_n$ is locally compact:
Then the spaces in (\ref{ctsncts})
coincide, and $\cT=\cO$ (see, e.g., \cite[Theorem~4.1]{Hir}).\\[2.5mm]
In this article, we discuss
more comprehensive classes of
topological groups
(and spaces),
for which the preceding difficulties
cannot occur.
Although the topological groups considered
need not be locally compact,
we obtain
results concerning
final group topologies,
colimits and direct limits
for such groups,
which involve compactness in less direct ways.
In particular,
we obtain results for topological groups
whose underlying topological spaces are
$k_\omega$ or locally~$k_\omega$.\\[3mm]
%
%
Recall that a topological space is called a
\emph{$k_\omega$-space} if it is the direct limit of an ascending
sequence of compact (Hausdorff) subspaces. The
class of $k_\omega$-spaces is quite
general and comprises, for example,
all countable CW-complexes~\cite[\S5]{Ord} and
countable direct limits of
$\sigma$-compact locally compact groups (cf.\ \cite{DIR},
\cite{TSH}); yet it is
restrictive enough to allow a systematic treatment. 
In the
theory of topological groups,
$k_\omega$-spaces have attracted considerable interest
because
free topological groups over compact (or
$k_\omega$-) spaces are $k_\omega$-spaces (see
\cite[Theorem~4]{Grv}, \cite[Corollary~1]{MMO} and
\cite[Theorem~5.2]{Or3}), and these are essentially the only
examples of free topological groups whose topology is well
accessible. The situation is similar for
Hausdorff quotient groups of such free groups
(which are $k_\omega$-spaces as well).
These subsume
important categorical constructs, like free products of two
$k_\omega$-groups with (or without) amalgamation
(\cite{Kat}--\cite{KM2}, \cite{KhM2}).\\[3mm]
The class of locally $k_\omega$ spaces
(introduced below) provides
a common roof for the classes of locally compact spaces on the one hand,
and $k_\omega$-spaces on the other.\\[3mm]
Our motivation to study $k_\omega$-spaces stems
from two specific sources of examples:
On the one hand, these are
examples from (infinite-dimensional) topological geometry and
geometric group theory, more precisely Moufang topological
twin buildings (cf.\ \cite{Hartnick}) and real and complex Kac-Moody
groups including their non-split forms and certain twisted
variants (cf.\ \cite{Kac/Peterson:1985}, \cite{TitsFunctor}
and \cite{Remy}, respectively).
On the other hand, $k_\omega$-groups
arise most naturally in the duality theory of
topological abelian groups, because dual groups of
metrizable abelian groups are~$k_\omega$
\cite[Corollary~4.7]{Aus}.\\[3mm]
The logical structure of the article guarantees
that readers interested
in only one of the two fields
of applications can skip the discussion
of the other.\\[4mm]
\noindent
{\bf Local compactness of colimits.}
After a preparatory section,
in Section~\ref{seclcp}
we study situations ensuring that a final
group topology is locally compact
(Proposition~\ref{praeprop}).
For instance,
we consider a diagram of countably many $\sigma$-compact
locally compact groups
and its colimit $G$
in the category of abstract groups.
If $G$
admits a locally compact group topology making
the limit maps continuous, then this
topology makes $G$ the colimit
in the category of topological groups
(Corollary~\ref{propo1}).
We remark that results concerning local compactness of
colimits for special diagrams are available
in the literature, notably for free products
$A*B$ and free products $A*_CB$ with amalgamation.
In these cases, local compactness
usually rules out interesting situations.%
\footnote{For example, for $A,B$ non-trivial,
local compactness of $A*B$
necessitates that $A$ and $B$ are discrete
(and likewise for
infinite free products,
see \cite[Corollary~1]{MaN}).
Furthermore, if $C$ is a proper subgroup
of both $A$ and~$B$,
then local compactness of $A*_CB$ forces
that $C$ is open in
$A*_CB$ (see \cite[Corollary~3]{AlT},
or \cite[Theorem~3]{KhM2}
for the case where~$C$ is normal
in~$A$ and $B$).
It is also known that an amalgamated product
$A*_CB$, with $C$ a proper subgroup
of~$A$ and~$B$, cannot be compact~\cite{AlC}.}
The picture changes if one considers more complicated
diagrams of topological groups
and their colimits (e.g.,
colimits of amalgams of more than two factors).
Then a rich supply of non-trivial
examples becomes available,
as we hope to
illustrate
by results concerning compact Lie groups
in Section~\ref{secphan}.\\[3mm]
We mention that
a typical feature of our
amalgamation results
becomes visible
at this point:
While topological group
theoreticians mainly studied
simple types of amalgams
of quite general groups,
our geometrically motivated
examples involve
complicated amalgams
of rather special groups.\\[4mm]
{\bf Lie groups determined by an amalgam of subgroups.}
Amalgams and their universal enveloping groups (i.e.,
their colimits in the category of groups)
play an important role in group theory.
The possibly most famous result in this area, the Curtis-Tits Theorem,
asserts that the universal version of an arbitrary Chevalley group
of (twisted) rank at least
three is the universal enveloping group of the amalgam consisting of the
groups $\langle X_{\pm \alpha} \rangle$ and
$\langle X_{\pm \alpha}, X_{\pm \beta} \rangle$,
where $\alpha$, $\beta$ run through a fundamental system of roots
corresponding to the Chevalley group and $X_\alpha$ denotes the
respective root subgroup
(cf.\ \cite{Curtis:1965}, \cite{Tits:1962},
also \cite[Theorem 2.9.3]{Gorenstein/Lyons/Solomon:1998}).
A similar result has been proved for compact Lie groups
(cf.\ \cite{Borovoi:1984}), and for complex Kac-Moody groups and their
unitary forms (cf.\ \cite{Kac/Peterson:1985}), as abstract groups.\\[2.5mm]
A lemma by Tits (cf.\ \cite{Tits:1986};
restated as Lemma \ref{tits}
in the present article) ties the\linebreak
above amalgamation results to the theory of simplicial complexes and,
in particular, to the theory of Tits
buildings and twin buildings.
For compact Lie groups and unitary
forms of complex Kac-Moody groups, this
relies on the fact that a Tits building of rank
at least three is simply connected
(cf.\ \cite{Tits:1974}).
For Chevalley groups and Kac-Moody groups,
it\linebreak
relies on the fact that the opposites
geometry of the corresponding twin
building is simply connected (cf.\ \cite{Muehlherr}).
Recently, a machinery was developed to classify the corresponding
amalgams, thus improving the above-mentioned results,
cf.\ \cite{Bennett/Shpectorov:2004}, \cite{D},
\cite{Gramlich:2004}, \cite{Gra}.\\[3mm]
The amalgamation results from \cite{Gra} treated compact Lie
groups as abstract groups. In Section~\ref{secphan}, we transfer
them into the setting of locally compact groups (and Lie groups),
using results from Section~\ref{seclcp}. All relevant mappings
identified earlier as isomorphisms of abstract groups turn out to
be topological isomorphisms (see Theorems~\ref{borovoi}
and~\ref{mainphan} and their respective proofs).\\[4mm]
{\bf Locally {\boldmath $k_\omega$} spaces.}
In Section~\ref{seckomeg},
we introduce locally $k_\omega$ spaces
(Hausdorff spaces each point of which has an open neighbourhood
which is a $k_\omega$-space).
These form a quite natural class of
topological spaces,
which subsumes all $k_\omega$-spaces
and all locally compact spaces
(and is therefore prone to unify results
known for the two subclasses).
Given an ascending sequence
$X_1\sub X_2\sub\cdots$
of locally $k_\omega$ spaces
with continuous inclusion maps
$X_n\to X_{n+1}$,
we show that the direct limit topology
on $X=\bigcup_{n\in \N}X_n$
is locally $k_\omega$
(Proposition~\ref{DLlkom}).
Furthermore, given another such sequence
$Y_1\sub Y_2\sub\cdots$, with direct limit~$Y$,
we show that the product topology
on $X\times Y$ makes $X\times Y$ the direct
limit topological space $\dl\, X_n\times Y_n$\vspace{-.3mm}
(Proposition~\ref{herecompat}).
The compatibility of direct limits and
direct products
was known previously for ascending sequences of locally
compact spaces (see, e.g., \cite[Theorem~4.1]{Hir} and
\cite[Proposition~3.3]{DIR}).
It is needed to get a grip on direct limits
of topological groups \cite{DIR}, \cite{Hir}
(as already indicated)
and direct limits of Lie groups
(see \cite{DIR}--\cite{DL3}, \cite{NRW} and \cite{NRW2}).\\[4mm]
{\bf Locally {\boldmath $k_\omega$} groups
and their direct limits.}
In Section~\ref{lkomegagp}, we
show that a topological group
is locally $k_\omega$
if and only if it has an open
subgroup which is a $k_\omega$-group
(Proposition~\ref{insplyd}).
We also prove that,
for each ascending sequence
$G_1\leq G_2\leq \cdots$ of locally~$k_\omega$ groups
with continuous inclusion maps,
the direct limit topology on
$G:=\bigcup_{n\in \N}G_n$
is locally $k_\omega$
and makes $G$ a topological
group (Proposition~\ref{injsequent}).
This result provides a\linebreak
common generalization
of known facts concerning direct limits
of locally compact groups (\cite[Corollary~3.4]{DIR},
\cite[Theorem~2.7]{TSH})
and abelian $k_\omega$-groups \cite[Proposition~2.1]{ATC}.
We also formulate a condition
ensuring that a final group
topology is $k_\omega$ (Proposition~\ref{whenfinkom}).
This entails the existence of certain colimits
of $k_\omega$-groups (Corollary~\ref{colimkom}).
For closely related results
concerning the group topology generated by a subspace~$X$
which is compact or $k_\omega$,
see \cite{MMO}, \cite{Mor} and \cite{Num}.\\[4mm]
{\bf Kac-Moody groups determined by an amalgam of subgroups.}
The topological amalgamation results discussed in
Section~\ref{secphan} have an infinite-dimensional analogue
which we describe in Section~\ref{SecKacMoody}.
On the algebraic level,
unitary forms of complex Lie groups are replaced by unitary forms
of complex Kac-Moody groups;
on the geometric level, spherical
buildings are replaced by twin buildings; and on the topological
level, (locally) compact spaces are replaced by $k_\omega$-spaces.
While the proofs become more technical, the results (cf. Theorem
\ref{TopCompletion}) are analogous to the
finite-dimensional case.\\[4mm]
{\bf The duality of locally {\boldmath $k_\omega$},
resp., almost metrizable abelian groups.}
Sections~\ref{secdualcat}
and~\ref{dirandpro}
are devoted to the duality theory
of (Hausdorff) topological abelian groups.
Given an abelian Hausdorff
group~$G$, let $G^*$ be its \emph{dual group},
i.e., the group of all continuous homomorphisms
from~$G$ to the circle group
$\T=\{z\in\C\colon |z|=1\}$, equipped with the
compact-open topology.
It is known that $G^*$
is a $k_\omega$-group
if~$G$ is metrizable
(see \cite[Corollary~4.7]{Aus}).
Conversely,
$G^*$ is metrizable (and complete)
for each topological abelian group~$G$
which is a $k_\omega$-space
by \cite[Propositions~2.8 and 4.11]{Aus}.
The duality of metrizable groups
and $k_\omega$-groups
was explored in~\cite{ATC},
in particular the duality between
countable direct limits of $k_\omega$-groups
and countable projective limits
of metrizable groups.
Our goal is to generalize
the results from~\cite{ATC} to the larger classes of
almost metrizable abelian groups
and locally~$k_\omega$ abelian groups.
Also these are in duality:
Exploiting the fact that
a topological abelian group~$G$ is
almost metrizable
(as in \cite{Pas} or \cite[Definition~1.22]{Aus})
if and only if it has a
compact subgroup
$K$ such that $G/K$ is metrizable
(see \cite[Proposition~2.20]{Aus} or \cite{Pas}),
we show in Section~\ref{secdualcat}
that $G^*$ is locally $k_\omega$ for
almost metrizable~$G$,
while $G^*$ is almost metrizable
whenever $G$ is locally $k_\omega$
(Proposition~\ref{dualcats}).\vfill\pagebreak

\noindent
{\bf Pontryagin duality for direct and projective limits.} In
Section~\ref{dirandpro}, we study the dual groups of countable
direct limits of locally $k_\omega$ abelian groups, and the dual
groups of countable projective limits
of almost metrizable abelian groups. Since locally~$k_\omega$ groups
subsume both locally compact abelian groups and abelian
$k_\omega$-groups, our results provide a common generalization of
Kaplan's classical results concerning the duality of countable
direct (and projective) limits of locally compact abelian groups
(\cite{Kp1}, \cite{Kp2}) and the recent studies from~\cite{ATC}
just mentioned. In the special cases of complete metrizable groups
and $k_\omega$-groups treated in~\cite{ATC}, we are able to loosen
some of the hypotheses of \cite{ATC} and thus\linebreak
to strengthen their results.
Cf.\ \cite{AT2} and~\cite{BaB}
for complementary
studies of the continuous duality of direct and
projective limits of convergence
groups. \\[4mm]
{\bf Acknowledgements.} The authors thank Andreas Mars (Darmstadt)
and an
anonymous referee for comments on a preliminary version of the
article.
The manuscript was
partially written while the second author held a substitute
professorship at the
TU Braunschweig. He would like to express his gratitude
to the members of the Institute of Computational Mathematics
of the TU Braunschweig
for their hospitality.
The first and second author acknowledge support
by the German Research Foundation (DFG)
(projects: GL 357/2-1 and GL 357/5-1,
respectively,
GR 2077/5-1 and GR 2077/7-1).
During the revision of this article, the third author was
partially
supported by Swiss National Science Foundation (SNF),
grant PP002-102765.
\section{Notation and terminology}\label{secprel}
Part of our terminology
has already been described in the introduction.
We now fix further terminology and notation.
In the following,
$\G$ denotes the category
of groups and homomorphisms.
All compact (or locally compact) spaces are assumed
Hausdorff.
The category of all (not necessarily
Hausdorff) topological groups
and continuous homomorphisms
will be denoted by
$\TG$,
and
$\HTG$, $\LCG$ and $\LIE$
are its full subcategories
of Hausdorff topological groups,
locally compact topological
groups and finite-dimensional
real Lie groups,
respectively. Surjective, open, continuous
homomorphisms
between topological groups shall be referred
to as \emph{quotient morphisms}.
Finally, we shall use
the term ``locally convex space''
as an abbreviation for ``locally convex topological vector space.''
\\[2.5mm]
{\bf Diagrams and colimits.}
Recall that a category
$\I$ is called {\em small\/}
if its class $\ob(\I)$ of objects
is a set.
A {\em diagram\/} in
a category~$\A$ is a covariant functor
$\delta \colon \I\to \A$
from a small category~$\I$
to~$\A$.
A {\em cone\/} over a diagram $\delta \colon \I\to \A$
is a natural transformation $\phi\colon
\delta \stackrel{\cdot}{\to} \gamma $ to a ``constant''
diagram~$\gamma$, i.e.,
to a diagram $\gamma \colon \I\to \A$
such that $G:=\gamma(i)=\gamma(j)$
for all $i,j\in \ob(\I)$
and $\gamma(\alpha)=\id_G$
for each morphism $\alpha$ in~$\I$.
We can think of a cone
as the object $G\in \ob(\A)$,
together with the family $(\phi(i))_{i\in \ob(\I)}$
of morphisms
$\phi(i)\colon \delta(i) \to G$,
such that
\[
\phi(j)\circ \delta(\alpha)\;=\; \phi(i)\qquad
\mbox{for all $\, i,j\in \ob(\I)$ and $\alpha\in \Mor(i,j)$.}
\]
A cone $(G,(\phi_i)_{i\in \ob(\I)})$
is called a {\em colimit\/} of $\delta$
if, for each cone $(H,(\psi_i)_{i\in \ob(\I)})$,
there is a unique
morphism $\psi\colon G\to H$ such that
$\psi\circ \phi_i=\psi_i$ for all $i\in \ob(\I)$.
If it exists,
a colimit is unique
up to natural isomorphism.
We expect that the reader is familiar
with the most basic types
of colimits and corresponding diagrams,
like direct limits (corresponding to the case
where $\I$ is a directed set),
free products of groups, and amalgamated
products $A*_CB$.\\[3mm]
{\bf Colimits of topological groups.}
Standard arguments
show that colimits
exist for any diagram
in $\G$, $\TG$ and $\HTG$
(as is well known).
Given a diagram
$\delta\colon \I\to\TG$
of topological groups,
with colimit $(G,(\lambda_i)_{i\in \ob(\I)})$
in the category $\G$ of abstract groups,
there is a finest group topology~$\cO$
on~$G$ making each $\lambda_i$ continuous;
then $((G,\cO),(\lambda_i)_{i\in \ob(\I)})$
is the colimit of $\delta$ in~$\TG$.
If $\delta$ is a diagram in $\HTG$ here,
then $H :=G/\wb{\{1\}}$
is a Hausdorff group
and $(H, (q\circ \lambda_i)_{i\in \ob(\I)})$
is readily seen to be
the colimit of $\delta$
in $\HTG$,
where $q\colon G\to G/\wb{\{1\}}=H$
is the canonical quotient morphism.
Frequently, the final group
topology on~$G$ is Hausdorff,
in which case the $\HTG$-colimit
is obtained by topologizing the colimit
of abstract groups.
For example, this situation occurs
in the case of free products
of Hausdorff topological groups \cite{Gv2}
and amalgamated products $A*_CB$
with $C$ a closed central subgroup of~$A$
and~$B$ (see \cite{KhM});
if $A$ and $B$ are $k_\omega$-groups,
it suffices that $C$ is a closed
normal subgroup~\cite{KaM}
or compact~\cite{KM1}.
This additional information is not available
through general category-theoretical arguments,
but requires a solid amount of work.\\[3mm]
{\bf Amalgams and their enveloping groups.}
In the article \cite{Gra} to be extended
in Section~\ref{secphan},
amalgams of groups play a central role,
which were defined there as families $(H_i)_{i\in J}$
of groups satisfying certain axioms,
ensuring in particular
that $H_i\cap H_j$
is a subgroup of $H_i$ and~$H_j$,
for all $i,j\in J$
(see \cite[Definition~2.5]{Gra}).
In the present paper, we
prefer a (wider) definition of amalgams
in category-theoretical terms:\\[2mm]
A diagram $\delta\colon \I\to \G$ of groups
is called an \emph{amalgam}
if $\I$ is the category associated with some partially
ordered set $(J,\leq)$ and $\delta(\alpha)$
is a monomorphism of groups
for each morphism~$\alpha$ in~$\I$.
The first condition means that
$\ob(\I)=J$,
and furthermore for $a,b\in J$ there exists one (and only one)
morphism $a\to b$ if and only if $a \leq b$. (Sometimes we prefer
to reverse the partial order on $J$ to the effect that there exists
a morphism $a \to b$ if and only if $a \geq b$).
Following the terminology in~\cite{Gramlich},
a cone $(G,(\phi_i)_{i\in \ob(\I)})$ over $\delta$ in~$\G$
will also be called an \emph{enveloping group} of the amalgam
and its colimit
a \emph{universal enveloping group}.
Analogous terminology will be used for amalgams of topological
groups. We would like to point out that in \cite{Goldschmidt:1980}, 
\cite{Gra} and other references, the terminology ``completion'' 
and ``universal completion'' has been used. This terminology, however, 
is unfortunate, because an amalgam of topological groups might 
have a universal enveloping group which is not a complete topological
group (in its left uniform structure).
As it may sometimes be useful in this case
to pass to the completed topological group, 
the clash of terminology might lead to confusion. \\[2mm]
Of particular importance for us are amalgams
$\delta\colon \I\to \G$,
where $J$ is a set and
$\I$ the small category
whose set of objects
is $\ob(\I)={ J \choose 1} \cup { J \choose 2}$,
and such that, besides identity morphisms,
there is exactly one morphism
$\{i\} \to \{i,j\}$ and $\{j\} \to \{i,j\}$
for all distinct $i,j\in J$.
We can think of such an amalgam as a family
of groups $(G_i)_{i \in \ob(\I)}$
(where $G_{\{i\}}:=\delta(\{i\})$,
$G_{\{i,j\}}:=\delta(\{i,j\})$)
together with monomorphisms of groups
$G_{\{i\}} \to G_{\{i,j\}}$, resp.,
$G_{\{j\}} \to G_{\{i,j\}}$
(cf.\ also \cite[Definition 2.15]{Gramlich/Maldeghem:2006}
for this type of amalgam).\\[3mm]
{\bf Basic facts concerning direct limits of topological spaces.}
For basic facts concerning direct limits
of topological spaces and topological groups,
the reader is referred to
\cite{DIR}, \cite{Han}, \cite{Hir}
and \cite{TSH}
(where also many of the pitfalls and subtleties
of the topic are described).
In particular, we shall frequently
use that the set underlying
a direct limit of groups is the
corresponding direct limit in the category
of sets, and that the direct limit
of a direct system of topological spaces
in the category of (not necessarily
Hausdorff) topological spaces
is its direct limit in the category
of sets, equipped with the final topology with respect
to the limit maps.
For later reference,
we compile further simple facts in a lemma:
%
%
\begin{la}\label{baseDL}
Let $((X_i)_{i\in I}, (f_{ij})_{i\geq j})$
be a direct system
in the category of topological spaces
and continuous maps,
with direct limit $(X,(f_i)_{i\in I})$.
\begin{itemize}
\item[\rm (a)]
If $U_i\sub X_i$ is open and $f_{ij}(U_j)\sub U_i$
for $i\geq j$,
then $\bigcup_{i\in I}f_i(U_i)$
is open in~$X$, and $U=\dl\,U_i$\vspace{-1mm}
holds for the induced topology.
\item[\rm (b)]
If $Y\sub X$ is closed or open, then
$Y=\dl \,f_i^{-1}(Y)$\vspace{-1mm}
as a topological space.
\item[\rm (c)]
If each $f_{ij}$ is injective,
then each $f_i$ is injective.
If each $f_{ij}$ is a topological
embedding, then each $f_i$
is a topological embedding.
\item[\rm (d)]
If $I$ is countable, each~$X_i$ is~$T_1$
and each $f_i$ is injective, then each compact subset~$K$
of~$X$ is contained in $f_i(X_i)$
for some $i\in I$.
If, furthermore, each $f_i$ is a
topological embedding, then
$K=f_i(L)$ for some $i\in I$
and some compact subset $L\sub X_i$.
\end{itemize}
\end{la}
\begin{proof}
(a) is quite obvious (cf.\ \cite[Lemma~3.1]{DIR}).
Part\,(b) follows from the definition
of final topologies.
The injectivity in (c) is clear
from the construction of~$X$
(see, e.g., \cite[\S2]{DIR}).
For the second assertion,
see \cite[Lemma~A.5]{NRW2}.
The first part of (d)
is \cite[Lemma~1.7\,(d)]{DL2},
and its second part follows
from\,(c) (see also \cite[Lemma~2.4]{Han}).
\end{proof}
If all bonding maps $f_{ij}\colon X_j\to X_i$
are topological embeddings,
then a direct system of topological spaces
(or topological groups) is called
\emph{strict}.
\section{Local compactness and final group topologies}\label{seclcp}
%
%
We describe a condition ensuring
that a final group topology is
locally compact.
A simple version
of the Open Mapping Theorem
will be used.
\begin{la}\label{openmap}
Let $f\colon G\to H$
be a surjective, continuous
homomorphism between Hausdorff
topological groups.
If $G$ is $\sigma$-compact
and $H$ is a Baire space,
then $f$ is open and
$H$ is locally compact.
\end{la}
\begin{proof}
By hypothesis,
$G=\bigcup_{n\in \N}K_n$
for certain compact sets $K_n\sub G$
and thus $H=\bigcup_{n\in \N} f(K_n)$
with $f(K_n)$ compact.
Since~$H$ is a Baire space,
$f(K_n)$ has non-empty interior
for some $n\in \N$, and thus
$H$ is locally compact.
Like any continuous surjection
between compact Hausdorff spaces,
$f|_{K_n} \colon K_n\to f(K_n)$
is a quotient map.
Let $q\colon G\to G/\ker(f)$
be the quotient homomorphism
and $\phi\colon G/\ker(f)\to H$ be the
bijective continuous homomorphism
induced by~$f$.
Then $\phi^{-1}\circ f|_{K_n}=q|_{K_n}$
is continuous,
whence $\phi^{-1}|_{f(K_n)}$ is continuous,
$\phi^{-1}$ is a continuous
homomorphism, and $\phi$ is a topological isomorphism.
\end{proof}
%
%
\begin{prop}\label{praeprop}
Let $G$ be a Hausdorff topological group
and $(f_i)_{i\in I}$ be a countable family
of continuous maps $f_i\colon X_i\to G$,
such that $X_i$ is $\sigma$-compact for each
$i\in I$ and $\bigcup_{i\in I}f_i(X_i)$
generates~$G$.
Then $G$ is locally compact if and only if~$G$
is a Baire space.
In this case,
the given locally compact topology on~$G$
is the final group topology with respect
to $(f_i)_{i\in I}$.
\end{prop}
\begin{proof}
Assume that $G$ is Baire
(which follows from local compactness).
The final group topology $\cT$
being finer than the given topology~$\cO$,
it is Hausdorff.
Since $(G,\cT)$
is generated by
its $\sigma$-compact
subset $\bigcup_{i\in I}f_i(X_i)$,
the topological group
$(G,\cT)$ is $\sigma$-compact.
As $(G,\cO)$ is Baire,
Lemma~\ref{openmap}
shows that
the continuous surjective
homomorphism
$(G,\cT)\to (G,\cO)$, $g\mto g$
is an isomorphism of topological
groups and $(G,\cO)$
is locally compact.
\end{proof}
The following special case
(applied to amalgams of finitely
many Lie groups)  
will be used presently to
transfer the purely algebraic considerations from~\cite{Gra}
into a topological (and Lie theoretic) context.\\[3mm]
The setting is as follows:
$\delta \colon \I\to \LCG$ is a diagram
of $\sigma$-compact locally compact groups
$G_i:=\delta(i)$ for $i\in I:=\ob(\I)$
and continuous homomorphisms
$\phi_\alpha:=\delta(\alpha) \colon G_i\to G_j$
for $i,j\in I$ and $\alpha\in \Mor(i,j)$,
such that $I$ is countable.
Furthermore, $(G,(\lambda_i)_{i\in I})$
is a colimit of the diagram~$\delta$
in the category of abstract groups,
with homomorphisms
$\lambda_i\colon G_i\to G$.
%
%
\begin{cor}\label{propo1}
If there exists
a locally compact Hausdorff group
topology $\cO$ on~$G$
making $\lambda_i\colon G_i\to (G,\cO)$
continuous for each $i\in I$, then
$((G,\cO),(\lambda_i)_{i\in I})$
is a colimit of~$\delta$
in the category
of topological groups,
in the category of Hausdorff groups,
and in the category
of locally compact groups.
In particular, there is only one
such locally compact group topology on~$G$.
If each $G_i$ is a $\sigma$-compact
Lie group and $(G,\cO)$
is a Lie group,
then $((G,\cO),(\lambda_i)_{i\in I})$ also is a
colimit of $\delta$ in the category $\LIE$
of Lie groups.\Punkt
\end{cor}
\section{Lie groups determined
by rank two subgroups}\label{secphan}
%
We are now in the position to
formulate a topological
variant of the main results of~\cite{Borovoi:1984}.\\[2.5mm]
Let $G$ be a simply connected compact semisimple Lie group,
let $T$ be a maximal torus of $G$, let
$\Sigma = \Sigma(G_{\mathbb{C}},T_{\mathbb{C}})$
be its root system, and let $\Pi$ be a system of fundamental roots
of $\Sigma$.
To each root $\alpha \in \Pi$ corresponds some compact semisimple
Lie group $K_\alpha \leq G$ of rank one such that $T$ normalizes
$K_\alpha$. For simple roots $\alpha$, $\beta$, we denote by
$K_{\alpha\beta}$ the group generated by the groups $K_\alpha$ and
$K_\beta$, and by $\Sigma_{\alpha\beta}$ its root system relative
to the torus $T_{\alpha\beta} = T \cap K_{\alpha\beta}$.
The group $K_{\alpha\beta}$ is a compact semisimple Lie group
of rank two and $\lb \alpha, \beta \rb$ is a fundamental system
of $\Sigma_{\alpha\beta}$.
The pair $(K_\alpha,K_\beta)$ is called a \emph{standard pair}
of $K_{\alpha\beta}$, as is any image of $(K_\alpha,K_\beta)$
under a homomorphism from $K_{\alpha\beta}$ onto a central
quotient of $K_{\alpha\beta}$. Standard pairs are conjugate.
In fact, maximal tori are conjugate (cf.\
\cite[Theorem~6.25]{Hofmann/Morris:1998}),
and, if $\alpha, \beta \in \Pi$ and $\alpha_1, \beta_1 \in \Sigma$
have the same lengths and the same angle, there exists an
element~$w$ of the Weyl group with $w(\alpha_1) = \alpha$
and $w(\beta_1) = \beta$ (cf.\ \cite{Bourbaki:2002}).\\[2.5mm]
Then the following topological version of Borovoi's theorem
\cite{Borovoi:1984} holds:
\begin{thm} \label{borovoi}
Let $G$ be a simply connected compact semisimple Lie group with
Lie group topology $\cO$, let $T$ be a maximal torus of $G$, let
$\Sigma = \Sigma(G_{\mathbb{C}},T_{\mathbb{C}})$ be its root system,
and let $\Pi$ be a system of fundamental roots of $\Sigma$.
Let $\I$ be a small category with objects
${\Pi\choose 1} \cup {\Pi \choose 2}$ and morphisms
$\{\alpha\} \to \{\alpha,\beta\}$,
for all $\alpha, \beta \in \Pi$, and let $\delta \colon \I \to \LCG$
be a diagram with $\delta(\{\alpha\}) = K_\alpha$,
$\delta(\{\alpha,\beta\}) = K_{\alpha\beta}$, and
$\delta(\{\alpha\} \to \{\alpha,\beta\})
= (K_\alpha \hookrightarrow K_{\alpha\beta})$.\\[2.5mm]
Then $((G, \cO), (\iota_i)_{i \in {\Pi\choose 1} \cup {\Pi \choose 2}})$,
where
$\iota_i$ is the natural inclusion map, is a colimit of $\delta$
in the category $\LIE$ of Lie groups.
\end{thm}
Before we prove the theorem, we remind the reader of the notion
of a Tits building and restate Tits' Lemma.
\begin{defn}
Following \cite[Theorem 5.2]{Tits:1974}, see also
\cite[Section 2]{Garland:1973} and \cite[Section 6.8]{Benson:1991},
the {\em Tits building $\Delta(G,\F)$ of $G$ over $\F$},
where $G$ is a reductive algebraic group defined over a field $\F$,
consists of the simplicial complex whose simplices are indexed by
the $\F$-parabolic subgroups of $G$ ordered by the reversed
inclusion relation on the parabolic subgroups.
\end{defn}
\begin{la}[Tits' Lemma; cf.\ \cite{Tits:1986}] \label{tits}
Let $A$ be a poset, let $G$ be a group acting on $A$, let
$F$ be a fundamental domain for the action of $G$ on $A$,
i.e., {\rm (i)} $a \in F$ and $b \leq a$ implies $b \in F$,
{\rm (ii)} $A = G.F$, {\rm (iii)} $G.a \cap F = \{ a \}$ for all
$a \in F$, and let $\I$ be a small category with objects $F$
and morphisms $b \leftarrow a$ for all $b \leq a \in F$.
Then the poset $A$ is simply connected if and only if $G$
is the colimit of the diagram $\delta \colon \I \to \G$,
where $\delta(a) = G_a$ and
$\delta(b \leftarrow a) = (G_b \hookleftarrow G_a)$.
\end{la}
Here a poset $A$ is called simply connected
if the associated simplicial complex $|A|$ -- with the non-empty
finite chains of $A$ as simplices -- is simply connected.\\[2.5mm]
\begin{Proof}{of Theorem \ref{borovoi}}
For Lie rank $|\Pi|$ of $G$ at most two there is nothing to show,
so we can assume that $|\Pi| \geq 3$.
Let $G_\C$ be the complexification of $G$ and let $\Delta$
be the Tits building of $G_\C$ (over $\C$),
which is a $(|\Pi|-1)$-dimensional simplicial complex.
Let $B$ be a Borel subgroup, i.e., a minimal parabolic subgroup,
of $G_\C$, containing the maximal torus $T$.
The $(|\Pi|-1)$-simplex $\sigma$ of $\Delta$ corresponding to $B$
is a fundamental domain for the conjugation action of $G_\C$
on $\Delta$ (see \cite[Theorem~3.2.6]{Tits:1974},
also \cite[Section 1]{Garland:1973}, \cite[Chapter 6]{Springer:1998}).
By the Iwasawa decomposition (see \cite[Theorem VI.5.1]{Helgason:1978}
or \cite[Theorem III.6.32]{Hilgert/Neeb:1991}),
we have $G_\C = GB$, so $\sigma$ is also a fundamental domain
for the conjugation action of $G$ on $\Delta$.
The stabilizers in $G$ of the 
sub-simplices codimension one and codimension two of $\sigma$  are exactly
the groups~$K_\alpha T$ and~$K_{\alpha\beta}T$.
By the simple connectedness of Tits buildings of rank
at least three
(cf.\ \cite[Theorem IV.5.2]{Brown:1989}, \cite[Theorem 13.32]{Tits:1974};
also \cite[Theorem 2.2]{Garland:1973}) plus an induction on~$|\Pi|$
applied to Lemma \ref{tits}, the group~$G$ is a colimit | in
the category of abstract groups | of the diagram
$\delta' \colon \I \to \LCG$
with $\delta'(\{\alpha\}) = K_\alpha T$,
$\delta'(\{\alpha,\beta\}) = K_{\alpha\beta}T$, and
$\delta'(\{\alpha\} \to \{\alpha,\beta\})
= (K_\alpha T \hookrightarrow K_{\alpha\beta} T)$.
By \cite[Lemma~29.3]{Gorenstein/Lyons/Solomon:1995}
(or by a reduction argument as in the proof of
\cite[Theorem 2]{Gramlich/Hoffman/Shpectorov:2003}
or in the proof of \cite[Theorem 4.3.6]{Gramlich:2004}),
the torus~$T$ can be reconstructed from the rank two tori
$T_{\alpha\beta}$, $\alpha, \beta \in \Pi$, and so the group~$G$
in fact is a colimit | in the category of abstract groups | of
the diagram~$\delta$.
Application of Corollary \ref{propo1} to $\delta$ yields the claim.
\end{Proof}
\section{Basic facts on {\boldmath $k_\omega$}-spaces
and locally {\boldmath $k_\omega$} spaces}\label{seckomeg}
%
%
%
%
In this section,
we introduce locally $k_\omega$ spaces
and compile various useful properties
of $k_\omega$-spaces and locally~$k_\omega$
spaces, for later use.
In particular, we discuss direct limits
in the category of locally~$k_\omega$ spaces.
%
%
\begin{defn}\label{defkomeg}
A Hausdorff topological
space~$X$ is a \emph{$k_\omega$-space}
if there exists
an ascending sequence of compact
subsets $K_1\sub K_2\sub \cdots \sub X$
such that $X=\bigcup_{n\in \N}K_n$
and $U \sub X$ is open
if and only if $U \cap K_n$ is open
in~$K_n$ for each $n\in \N$
(i.e., $X=\dl\, K_n$\vspace{-.3mm} as a topological space).
Then $(K_n)_{n\in \N}$
is called a \emph{$k_\omega$-sequence}
for~$X$.
We say that a Hausdorff topological space~$X$ is
\emph{locally $k_\omega$}
if each point has an open neighbourhood
which is a $k_\omega$-space in the induced topology.
\end{defn}
Thus every $k_\omega$-space is locally $k_\omega$.
It is clear that every locally $k_\omega$-space~$X$ is a $k$-space
(i.e., $X$ is Hausdorff and a subset $U\sub X$ is open
if and only if $U\cap K$ is open in~$K$ for each
compact set $K\sub X$).
%
%
\begin{prop}\label{basekomeg}
\begin{itemize}
\item[\rm (a)]
$\sigma$-compact locally compact spaces
are $k_\omega$.
\item[\rm (b)]
Closed subsets of $k_\omega$-spaces are
$k_\omega$ in the induced topology.
\item[\rm (c)]
Finite products of $k_\omega$-spaces are $k_\omega$
in the product topology.
\item[\rm (d)]
Hausdorff quotients of $k_\omega$-spaces
are $k_\omega$.
\item[\rm (e)]
Countable disjoint unions
of $k_\omega$-spaces are $k_\omega$.
\item[\rm (f)]
Every locally compact space is
locally $k_\omega$.
\item[\rm (g)]
Open subsets of locally $k_\omega$ spaces
are locally $k_\omega$.
\item[\rm (h)]
Closed subsets of locally $k_\omega$ spaces are
locally $k_\omega$ in the induced topology.
\item[\rm (i)]
Finite products of locally $k_\omega$
spaces are locally $k_\omega$ in the product topology.
\item[\rm (j)]
Hausdorff quotients of locally $k_\omega$ spaces
under open quotient maps are locally
$k_\omega$.
\end{itemize}
\end{prop}
\begin{proof}
(a)--(e) are well known:
see, e.g., \cite{FaT}
and the references therein.\footnote{For the reader's convenience,
we mention: In (a), any
exhaustion $K_1\sub K_2\sub\cdots$
of the space by compact sets
with $K_n$ in the interior
of $K_{n+1}$ provides a $k_\omega$-sequence.
Part\,(b) is a special case of Lemma~\ref{baseDL}\,(b).
Part\,(c) follows from Proposition~\ref{herecompat} below.
Part\,(d) follows from the transitivity of final
topologies.
Part\,(e): If $X=\coprod_{j\in \N}X_j$,
where each $X_j$ is a $k_\omega$-space with a
$k_\omega$-sequence $(K^{(j)}_n)_{n\in \N}$,
then clearly $K_n:=\bigcup_{j\leq n}K^{(j)}_n$
defines a $k_\omega$-sequence for~$X$.}
Part\,(h) follows from (b),
(i) from (c), and (j) from (d).

(f) If $X$ is locally compact and $x\in X$,
choose compact subsets $K_1, K_2,\ldots$
such that $x\in K_1$ and $K_{n+1}$
contains $K_n$ in its interior.
Then $U:=\bigcup_{n\in \N}K_n$
is an open neighbourhood of~$x$ in~$X$ and~$U$
is a $k_\omega$-space
with $(K_n)_{n\in \N}$ as a $k_\omega$-sequence.

(g) It suffices to show that each open neighbourhood~$W$
of a point~$x$ in a $k_\omega$-space~$X$
contains an open neighbourhood~$U$ of~$x$
which is a $k_\omega$-space.
To see this, let $(K_n)_{n\in \N}$
be a $k_\omega$-sequence for~$X$,
with $x\in K_1$.
Then $K_1\cap W$ contains a compact
neighbourhood $C_1$ of~$x$ in~$K_1$.
Recursively, we find a compact
subsets $C_n\sub W\cap K_n$
such that $C_n$
is a neighbourhood of $C_{n-1}$ in
$K_n\cap W$, for each $n\in \N$.
Let $U_n$ be the interior
of $C_n$ in~$K_n$.
Then $U:=\bigcup_{n\in \N} U_n\sub W$
is open in~$X$ and
$U=\dl\, U_n$\vspace{-.3mm}
in the induced topology, by Lemma~\ref{baseDL}\,(a).
Since
$U_n\sub C_n\sub U_{n+1}$ for each $n\in \N$,
we see that $U=\dl\,U_n=\dl\,C_n$\vspace{-1mm}
is a $k_\omega$-space.
\end{proof}
Proposition~\ref{basekomeg}\,(j)
ensures that
every Hausdorff quotient group
of a locally $k_\omega$ group
is locally $k_\omega$.
%
%
\begin{la}\label{somewstrange}
Let $X$ be a locally $k_\omega$ space
and $Y\sub X$ be a $\sigma$-compact subset.
Then $Y$ has an open neighbourhood
$U$ in~$X$ which is a $k_\omega$-space.
\end{la}
\begin{proof}
Since $Y$ is $\sigma$-compact,
we find a sequence $(U_n)_{n\in \N}$
of open subsets of~$X$
such that $Y\sub \bigcup_{n\in \N}U_n=:U$
and $U_n$ is a $k_\omega$-space,
for each $n\in \N$.
Then $U$ is a Hausdorff quotient
of the $k_\omega$-space $\coprod_{n\in \N}U_n$
and therefore a $k_\omega$-space,
by Proposition~\ref{basekomeg}\,(d) and\,(e).
\end{proof}
Recall that a map $f\colon X\to Y$
between topological spaces is called
\emph{compact-covering} if for
each compact subset $K\sub Y$,
there exists a compact set $L\sub X$
such that $K\sub f(L)$.
%
\begin{la}\label{coco1}
Every quotient map $q\colon X\to Y$
between $k_\omega$-spaces
is compact-covering.
\end{la}
\begin{proof}
Let $K\sub Y$ be compact.
If $(K_n)_{n\in \N}$
is a $k_\omega$-sequence for~$X$,
then $(q(K_n))_{n\in \N}$
is a $k_\omega$-sequence for~$Y$
and thus $K\sub q(K_n)$ for some~$n$,
by Lemma~\ref{baseDL}\,(d).
\end{proof}
%
%
\begin{prop}\label{DLlkom}
Let $X_1\sub X_2\sub\cdots$
be a sequence of $k_\omega$-spaces $($resp.,
locally~$k_\omega$ spaces$)$ $X_n$,
such that the inclusion map
$X_n\to X_{n+1}$ is continuous,
for each $n\in \N$.
Then the final topology $\cO$ on
$X:=\bigcup_{n\in \N}X_n$
with respect to the inclusion maps
$X_n\to X$ is $k_\omega$
$($resp., locally $k_\omega)$,
and makes $X$ the direct limit
$X=\dl\, X_n$\vspace{-.5mm}
in each of the categories of topological spaces,
Hausdorff topological spaces,
and $k_\omega$-spaces $($resp.,
locally $k_\omega$ spaces$)$.
\end{prop}
\begin{proof}
If each $X_j$ is a $k_\omega$-space,
let $(K_n^{(j)})_{n\in \N}$
be a $k_\omega$-sequence for $X_j$.
After replacing $K_n^{(j)}$
with $\bigcup_{i=1}^j K_n^{(i)}$,
we may assume that $K_n^{(i)}\sub K_n^{(j)}$
for all $i,j\in \N$ such that $i\leq j$.
Define $K_n:=K^{(n)}_n$ for $n\in \N$.
Since each $K_n$ is compact
and the inclusion maps $K_n\to K_{n+1}$
are topological embeddings,
the final topology $\cT$
on $X$ with respect to the inclusion maps
$K_n\to X$ is Hausdorff
(see \cite[Proposition~3.6\,(a)]{DIR}).
Thus $(X,\cT)=\dl\, K_n$\vspace{-.3mm} is a $k_\omega$-space.
The inclusion maps $K_n\to X_n\to (X,\cO)$
being continuous, we deduce
that $\cO\sub \cT$.
To see the converse, let $U\in \cT$.
Then $U\cap K_n$ is open in $K_n$
for each $n\in \N$.
Given $n\in \N$ and $j\in \N$,
let $m:=\max\{n,j\}$.
Then $U\cap K_n^{(j)}=
(U\cap K_m) \cap K_n^{(j)}$
is open in $K_n^{(j)}$,
the inclusion map $K^{(j)}_n\to K_m$ being continuous.
Hence $U\cap X_j$ is open in~$X_j$
for each~$j$ and thus
$U\in \cO$. Hence $\cO=\cT$, as required.\\[3mm]
In the case of locally $k_\omega$ spaces~$X_n$,
let $x\in X$, say $x\in X_1$
(after passing to a cofinal subsequence).
Let $U_1\sub X_1$ be an open neighbourhood
of $x$ which is a $k_\omega$-space.
Using Lemma~\ref{somewstrange},
we find open subsets $U_n\sub X_n$
which are $k_\omega$-spaces
and such that $U_n\sub U_{n+1}$,
for each $n\in \N$.
Then $U:=\bigcup_{n\in \N}U_n$
is an open subset of~$X$
(see Lemma~\ref{baseDL}\,(a))
and $U=\dl\,U_n$ is a $k_\omega$-space,
by what has just been shown.
It only remains to show that $X$ is Hausdorff.
To this end, suppose that $x,y\in X$
are distinct elements.
We may assume that $x,y\in X_1$.
The set $\{x,y\}$ being compact,
it has an open neighbourhood
$U_1\sub X_1$ which is a $k_\omega$-space
(by Lemma~\ref{somewstrange}).
Choose $U_n\sub X_n$ as just described.
Then both $x$ and $y$
are contained in the open subset
$U=\bigcup_{n\in\N} U_n$ of~$X$ which
is a $k_\omega$-space and hence Hausdorff,
whence $x$ and $y$ can be separated by open
neighbourhoods in~$U$.
\end{proof}
\begin{rem}
Hausdorff quotients of $k_\omega$-spaces
being $k_\omega$, it easily follows from
Proposition~\ref{DLlkom}
that the direct limit Hausdorff
topological space of each
countable direct system
of $k_\omega$-spaces is
a $k_\omega$-space.
\end{rem}
To formulate a result
concerning the compatibility of
direct limits and direct products,
let $X_1\sub X_2\sub\cdots$
and $Y_1\sub Y_2\sub\cdots$
be ascending sequences of topological
spaces with continuous inclusion maps,
and $X:=\bigcup_{n\in \N} X_n=\dl\,X_n$,
$Y:=\bigcup_{n\in \N}Y_n$,\vspace{-.3mm}
equipped with the direct limit topology.
Write $\dl\, (X_n\times Y_n)$\vspace{-.3mm}
for $\bigcup_{n\in \N}(X_n\times Y_n)$,
equipped with the direct limit topology.
%
%
\begin{prop}\label{herecompat}
The natural map
\[
\beta\colon \dl\,(X_n\times Y_n)
\to \big(\dl\,X_n\big)\times \big(\dl\, Y_n\big)\,,
\quad (x,y)\mto (x,y)
\]
is a continuous bijection.
If each $X_n$ and each $Y_n$ is locally $k_\omega$,
then $\beta$ is a homeomorphism.
\end{prop}
\begin{proof}
The inclusion map
$X_n\times Y_n\to X\times Y$
is continuous for each $n\in \N$.
Therefore~$\beta$ is continuous.
If each $X_n$ and $Y_n$ is locally $k_\omega$,
let $(x,y)\in X\times Y$,
say $(x,y)\in X_1\times Y_1$.
We show that there exist open neighbourhoods
$U\sub X$ of~$x$ and $V\sub Y$ of~$y$ such
that $X\times Y$ and $\dl\,(X_n\times Y_n)$\vspace{-.3mm}
induce the same topology on $U\times V$.
To this end, let $U_1$
be an open neighbourhood of~$x$
in~$X_1$ which is $k_\omega$, and choose
open subsets $U_n\sub X_n$
such that $U_n\sub U_{n+1}$
for each $n$ and each $U_n$ is~$k_\omega$
(using Lemma~\ref{somewstrange}).
Likewise, choose open, $k_\omega$
subsets $V_n\sub Y_n$ such that $y\in V_1$ and $V_n\sub
V_{n+1}$.
Then $U:=\bigcup_{n\in \N}U_n$
and $V:=\bigcup_{n\in \N}V_n$
are $k_\omega$-spaces in the topology induced
by~$X$ and~$Y$, respectively
(as a consequence of Proposition~\ref{basekomeg}\,(d) and\,(e)).
As in the first half of the proof of
Proposition~\ref{DLlkom},
we find $k_\omega$-sequences
$(K_n^{(j)})_{n\in \N}$ for each $U_j$
such that $K_n^{(i)}\sub K_n^{(j)}$
for $i\leq j$, whence $K_n:=K_n^{(n)}$
defines a $k_\omega$-sequence for~$U$.
Choose $k_\omega$-sequences $(L_n^{(j)})_{n\in \N}$
for $Y_j$
analogously
and define $L_n:=L_n^{(n)}$.
By the compatibility of direct products
and countable direct limits of compact spaces
(see, e.g., \cite[Proposition~3.3]{DIR}),
we have
\[
U_j\times V_j\,=\, \big(\dl\,K^{(j)}_n\big)\times
\big(\dl\, L^{(n)}_j\big)\,=\,
\dl\,(K_n^{(j)}\times L_n^{(j)})\,,\vspace{-1mm}
\]
showing that $(K_n^{(j)}\times L_n^{(j)})_{n\in \N}$
is a $k_\omega$-sequence for
$U_j\times V_j$.
Hence $(K_n\times L_n)_{n\in \N}$
is a $k_\omega$-sequence for
the open subset $\dl\,(U_n\times V_n)$
of $\dl\,(X_n\times Y_n)$,\vspace{-.3mm}
equipped with the induced topology.
Then
\[
U\times V\, =\, (\dl\, K_n)\times(\dl\,L_n)
\, =\, \dl\, (K_n\times L_n)\, =\, \dl\,(U_n\times V_n)\, ,\vspace{-1mm}
\]
using
\cite[Proposition~3.3]{DIR} for the second equality.
Here $\dl\,(U_n\times V_n)$\vspace{-.3mm}
is the set $U\times V$, equipped with the topology induced
by $\dl\, (X_n\times Y_n)$ (see Lemma~\ref{baseDL}\,(a)).
Thus $X\times Y$ and $\dl\, (X_n\times Y_n)$\vspace{-,3mm}
induce the same topology on this set,
which completes the proof.
\end{proof}
With a view towards later applications to Kac-Moody groups,
we recall a simple fact from \cite[\S\,21]{FaT}.
%
It illustrates that metrizability is
not preserved under passage to direct limits.
\begin{prop}\label{Metrizable}
If a locally $k_\omega$ space $X$ is metrizable,
then it is locally compact.\,\Punkt
\end{prop}
%
%
%
%
%
%
%
%
\section{Locally {\boldmath $k_\omega$} groups
and their direct limits}\label{lkomegagp}
%
%
%
In this section, we introduce and discuss the notion
of a locally $k_\omega$ group.
We also prove the existence of
direct limits for ascending
sequences of locally~$k_\omega$ groups,
and related results.
\begin{defn}
A \emph{$k_\omega$-group} (resp.,
locally $k_\omega$ group)
is a topological group the underlying topological
space of which is a $k_\omega$-space
(resp., a space which is locally $k_\omega$).
\end{defn}
\begin{ex}
Among the main examples of infinite-dimensional
Lie groups are direct limits $G=\dl\,G_n$\vspace{-.3mm}
of ascending sequences $G_1\sub G_2\sub\cdots$
of finite-dimensional Lie groups,
such that the inclusion maps are continuous
homomorphisms (see \cite{DIR}, \cite{DL2}, \cite{NRW}).
As $G$ is equipped with the final
topology here, Proposition~\ref{DLlkom}
shows that the topological space underlying~$G$
is locally $k_\omega$.
Thus $G$ is a locally $k_\omega$ group.
The same conclusion holds for
countable direct limits of locally compact groups,
as considered in \cite{DIR} and \cite{TSH}.
If each $G_n$ is $\sigma$-compact,
then~$G$ is a $k_\omega$-group.
However,
if some $G_n$ fails to be $\sigma$-compact
and each inclusion map is a topological embedding,
then~$G$ is not a $k_\omega$-group.
\end{ex}
%
%
%
\begin{prop}\label{insplyd}
For a topological group $G$,
the following conditions are equivalent:
\begin{itemize}
\item[\rm (a)]
$G$ is a locally $k_\omega$ group;
\item[\rm (b)]
$G$ has an open subgroup $H\leq G$
which is a $k_\omega$-group;
\item[\rm (c)]
$G=\dl\,X_n$\vspace{-.3mm} as a topological space
for an ascending sequence $X_1\sub X_2\sub\cdots$
of closed, locally compact subsets~$X_n$ of~$G$,
equipped with the induced topology.
\end{itemize}
\end{prop}
\begin{proof}
(a)$\impl$(b):
Let $U_1$ be an open identity neighbourhood
of~$G$ which is a $k_\omega$-space.
Then $U_1U_1 \cup U_1^{-1}$
is a $\sigma$-compact subset of~$G$ and hence has
an open neighbourhood $U_2$
which is a $k_\omega$-space.
Proceeding in this way, we obtain
an ascending sequence $U_1\sub U_2\sub\cdots$
of open identity neighbourhoods $U_n$
which are $k_\omega$-spaces, and such that
$U_nU_n\cup U_n^{-1}\sub U_{n+1}$
for each $n\in \N$.
Then $U:=\bigcup_{n\in \N}U_n$
is a subgroup of~$G$.
Apparently,
$U$ is open in~$G$ and
$U=\dl\,U_n$\vspace{-.3mm},
which is a $k_\omega$-space by Proposition~\ref{DLlkom}.

(b)$\impl$(c): Let $H\sub G$ be an open subgroup
which is a $k_\omega$-group,
and $(K_n)_{n\in \N}$ be a $k_\omega$-sequence
for~$H$. Let $T\sub G$ be a transversal
for $G/H$. Then $X_n:=\bigcup_{g\in T}gK_n$
is a closed subset of~$G$,
since $(gK_n)_{g\in T}$ is a locally finite
family of closed sets (see \cite[\S8.5, Hilfssatz~2]{Sch}).
Furthermore, $X_n$ is locally compact,
because $gK_n=X_n\cap gH$
is open in $X_n$.
Then $G=\bigcup_{n\in \N} X_n$
and $X_1\sub X_2\sub\cdots$.
The final topology $\cT$ on~$G$
with respect to the inclusion maps
$X_n\to G$ is finer than the original topology~$\cO$.
If $U\in \cT$,
let $g\in T$.
Then $(U\cap gH) \cap gK_n=(U\cap X_n)\cap gK_n$ is open
in $gK_n$ for each $n\in \N$.
Since $(gK_n)_{n\in \N}$
is a $k_\omega$-sequence for $gH$,
this entails that
$U\cap gH$ is open in $gH$
and hence open in~$G$.
Hence $U=\bigcup_{g\in T}(U\cap gH)\in \cO$.
Thus $\cO=\cT$.

(c)$\impl$(a): Locally compact spaces
being locally $k_\omega$,
(c) entails that $G$ is a direct limit
of an ascending sequence of locally $k_\omega$ spaces
and hence locally $k_\omega$, by Proposition~\ref{DLlkom}.
\end{proof}
%
%
%
\begin{prop}\label{injsequent}
Let $G_1\sub G_2\sub \cdots$
be an ascending sequence of $k_\omega$-groups
$($resp., locally $k_\omega$ groups$)$
such that the inclusion maps $G_n\to G_{n+1}$
are continuous homomorphisms.
Equip $G:=\bigcup_{n\in \N}G_n$
with the unique group structure making
each inclusion map $G_n\to G$ a homomorphism.
Then the final topology with respect to the
inclusion maps $G_n\to G$
turns $G$ into a $k_\omega$-group $($resp., a locally $k_\omega$
group$)$ and makes $G$ the direct limit
$\dl\,G_n$\vspace{-.3mm}
in the category of topological spaces,
and in the category of topological groups.
\end{prop}
\begin{proof}
The direct limit property in the category
of topological spaces is clear,
and we know from Proposition~\ref{DLlkom}
that~$G$ is a $k_\omega$-space
(resp., a locally $k_\omega$ space).
Let $\mu_n\colon G_n\times G_n\to G_n$
and $\lambda_n\colon G_n\to G_n$
be the group multiplication
(resp., inversion) of~$G_n$,
and $\mu\colon G\times G\to G$
and $\lambda \colon G\to G$
be the group multiplication (resp., inversion)
of~$G$.
Then $\lambda=\dl\,\lambda_n$,
where each $\lambda_n$ is continuous, and hence
$\lambda$ is continuous.
Identifying $G\times G$
with $\dl\, (G_n\times G_n)$\vspace{-.3mm}
as a topological space
by means of Proposition~\ref{herecompat},
the group multiplication~$\mu$ becomes
the map $\dl\,\mu_n\colon \dl\, (G_n\times G_n)\to\dl\, G_n$,
which is continuous.
Hence $G$ is a topological group.
The remainder is clear.
\end{proof}
%
%
\begin{cor}\label{beyondinj}
Let $I$ be a countable directed set
and $\cS:=((G_i)_{i\in I}, (f_{ij})_{i\geq j})$
be a direct system of
$k_\omega$-groups
$($resp., locally $k_\omega$ groups$)$
and continuous homomorphisms.
Then
$\cS$ has a direct limit $(G,(f_i)_{i\in I})$
in the category of Hausdorff topological
groups.
The topology on~$G$ is the final topology
with respect to the family $(f_i)_{i\in I}$,
and the Hausdorff topological space underlying
$G$ is the direct limit of~$\cS$ in
the category of Hausdorff topological
spaces. Furthermore, $G$ is a $k_\omega$-group
$($resp., locally
$k_\omega$ group$)$
and hence is the direct limit
of $\cS$ in the category of $k_\omega$-groups
$($resp., locally $k_\omega$ groups$)$.
\end{cor}
\begin{proof}
We may assume that $I =(\N,\leq)$.
Let $N_m:=\wb{\bigcup_{n\geq m}\ker f_{n,m}}$,
$Q_m:=G_m/N_m$ and $q_m\colon G_m\to Q_m$
be the canonical quotient morphism.
By Proposition~\ref{basekomeg}\,(d)
(resp., (j)),
each $Q_m$ is a $k_\omega$-group
(resp., a locally $k_\omega$ group).
Let $g_{n,m}\colon Q_m\to Q_n$
be the continuous homomorphism
determined by $g_{n,m}\circ q_m=q_n\circ f_{n,m}$.
Assume that $(X,(h_n)_{n\in \N})$
is a cone of continuous maps $h_n\colon G_n\to X$
to a Hausdorff topological space.
Then each $h_n$ factors to a continuous
map $k_n\colon Q_n\to X$.
{}From the preceding, it is clear
that the direct limit
of $\cS$ in the category of
Hausdorff topological spaces
(resp., Hausdorff topological groups)
coincides with that of
$((Q_n)_{n\in \N},(g_{n,m})_{n\geq m})$.
As each $g_{n,m}$ is injective,
the assertions follow from
Proposition~\ref{injsequent}.
\end{proof}
\begin{rem}
If each of the groups
$G_i$ in Proposition~\ref{beyondinj}
is abelian, then also~$G$ is abelian,
whence $G=\dl\,G_i$\vspace{-.3mm}
in the category of abelian $k_\omega$-groups
(resp., abelian, locally $k_\omega$ groups).
Thus
\cite[Prop.\,2.1]{ATC}
is obtained as special case,
the proof of which given in \cite{ATC}
is incorrect.\footnote{The problem
is that the final group topology on
$S:=\bigoplus_{n\in \N}G_n$
with respect to the inclusion maps
$i_n\colon G_n\to S$
does not coincide with the final
topology with respect to these
maps.}
\end{rem}
Inspired by arguments from the proof
of \cite[Proposition~2.1]{ATC}, we record:
%
\begin{cor}\label{dsums}
Let $(G_n)_{n\in \N}$ be a sequence
of abelian $k_\omega$-groups
$($resp., abelian, locally $k_\omega$ groups$)$.
Then the box topology
on $G:=\bigoplus_{n\in \N}G_n$
makes~$G$
the direct limit $\dl\,\prod_{k=1}^nG_n$\vspace{-.7mm}
in the category of topological
spaces, and this topology is $k_\omega$
$($resp., locally $k_\omega)$.
\end{cor}
\begin{proof}
By Proposition~\ref{basekomeg}\,(c)
(resp., (i)), $\prod_{k=1}^nG_k$ is
a $k_\omega$-group (resp., locally $k_\omega$).
By Proposition~\ref{injsequent},
the direct limit topology
turns $\bigoplus_{n\in \N}G_n=\dl\,\prod_{k=1}^nG_k$
into a $k_\omega$-group (resp., locally $k_\omega$ group),
and it makes $\bigoplus_{n\in\N}G_n$
the direct limit
$\dl\, \prod_{k=1}^nG_k$\vspace{-.3mm}
in the category $\TG$ of topological
groups. But it is well known
that the topology making
$\bigoplus_{n\in\N}G_n$
the direct limit
$\dl\, \prod_{k=1}^nG_k$\vspace{-.4mm}
in~$\TG$ is the box topology
(see, e.g., \cite[Lemma~4.4]{DL3}).
\end{proof}
Cf.\ also \cite{CaD} for recent investigations
of topologies on direct sums.\\[3mm]
The following proposition is a variant
of \cite[Theorem~1]{MMO}, which corresponds
to the most essential case where the final group
topology with respect to a single map
$f\colon X\to G$
from a $k_\omega$-space~$X$ to~$G$ is Hausdorff
and makes $f$ a topological embedding
(cf.\ also \cite{Num}, which even subsumes
non-Hausdorff situations).
%
%
\begin{prop}\label{whenfinkom}
Let $G$ be a group
and $(f_i)_{i\in I}$ be a countable
family of maps $f_i\colon X_i\to G$
such that each $X_i$ is a $k_\omega$-space
and $\bigcup_{i\in I}f_i(X_i)$
generates~$G$.
If the final group topology on~$G$
with respect to the family $(f_i)_{i\in I}$
is Hausdorff, then it makes $G$
a $k_\omega$-group.
\end{prop}
\begin{proof}
The disjoint union
$X:=\coprod_{i\in I}X_i$
is a $k_\omega$-space
(by Proposition~\ref{basekomeg}\,(e))
and hence a regular topological space
(see \cite[Proposition~4.3\,(ii)]{Han}).
We can therefore form
the Markov free topological group~$F$
on~$X$, as in \cite{Mar}.
The final group topology on~$G$
turns $G$ into a Hausdorff quotient group
of~$F$. Since~$X$ is a $k_\omega$-space,
$F$ is a $k_\omega$-group \cite[Corollary~1\,(a)]{MMO}.
Hence also its
Hausdorff quotient group~$G$
is a $k_\omega$-group, by
Proposition~\ref{basekomeg}\,(d).
\end{proof}
\begin{rem}
Note that, in the situation
of Proposition~\ref{whenfinkom},
the final group topology on~$G$
is Hausdorff if there exists a
Hausdorff group topology on~$G$
making all maps $f_i$ continuous.
\end{rem}
We record a direct consequence concerning colimits
in the category
$\KOG$ of $k_\omega$-groups.
%
%
\begin{cor}\label{colimkom}
Let $\delta \colon \I\to \KOG$ be a diagram
of $k_\omega$-groups
$G_i:=\delta(i)$ for $i\in I:=\ob(\I)$
and continuous homomorphisms
$\phi_\alpha:=\delta(\alpha) \colon G_i\to G_j$
for $i,j\in I$ and $\alpha\in \Mor(i,j)$,
such that $I$ is countable.
Let
$(G,(\lambda_i)_{i\in I})$
be a colimit of $\delta$ in the category
of Hausdorff topological groups,
with homomorphisms
$\lambda_i\colon G_i\to G$.
Then $G$ is a $k_\omega$-group,
and $(G,(\lambda_i)_{i\in I})$
also is the colimit of~$\delta$
in the category of $k_\omega$-groups.\Punkt
\end{cor}
%
%
\begin{la}\label{co-co}
Quotient morphisms
between locally~$k_\omega$
groups are compact-covering.
\end{la}
\begin{proof}
Given a quotient morphism $q\colon G\to Q$,
let $U\sub G$ be an open subgroup
which is a $k_\omega$-group.
Then $q(U)$ is an open subgroup
of~$Q$ and~$k_\omega$,
by Lemma~\ref{basekomeg}\,(d).
Given a compact set $K\sub Q$,
there are $y_1,\ldots,y_n\in K$ such
that $K\sub \bigcup_{j=1}^n y_jq(U)$.
By Lemma~\ref{coco1},
$y_j^{-1}(K\cap y_jq(U))\sub q(L_j)$
for a compact set $L_j\sub U$ and thus
$K\sub q(L)$ with $L:=\bigcup_{j=1}^n x_jL_j$,
where $x_j\in q^{-1}(\{y_j\})$.
\end{proof}
%
%
\begin{la}\label{findco}
In the situation of Corollary~{\rm \ref{beyondinj}},
let $K\sub G$ be a compact subset.
Then $K\sub f_i(L)$ for some $i\in I$
and some compact subset $L\sub G_i$.
\end{la}
\begin{proof}
Since the quotient morphisms
$q_i\colon G_i\to Q_i$ in the proof
of Corollary~\ref{beyondinj}
are compact-covering by Lemma~\ref{co-co},
we may assume that each $f_{ij}$
(and hence each $f_i$) is injective.
It therefore suffices to consider
the situation of Proposition~\ref{injsequent}.
Using Lemma~\ref{somewstrange}
and the argument from the proof
of ``(a)$\impl$(b)'' in Proposition~\ref{insplyd},
for each $n\in \N$ we find an open subgroup $U_n\sub G_n$
which is a $k_\omega$-group,
such that $U_n\sub U_{n+1}$ for each $n\in \N$.
Then $U:=\bigcup_{n\in \N}U_n$
is an open subgroup of $G=\bigcup_{n\in \N}G_n$
(see Lemma~\ref{baseDL}\,(a)), and the induced
topology makes $U$ the direct limit
$U=\dl\,U_n$.\vspace{-.8mm}
By Proposition~\ref{DLlkom}, $U$ is a $k_\omega$-space.
Choose a $k_\omega$-sequence
$(K^{(j)}_n)_{n\in \N}$
for $U_j$, for each $j\in \N$,
such that $K_n^{(i)}\sub K_n^{(j)}$
whenever $i\leq j$.
Then $(K_n^{(n)})_{n\in \N}$
is a $k_\omega$-sequence for~$U$.
Given a compact subset $K\sub U$,
we have $K\sub K_n^{(n)}$
for some $n\in \N$ (Lemma~\ref{baseDL}\,(d)),
whence~$K$ is a compact subset of~$U_n$.
Since every compact subset
of~$G$ is contained in finitely
many translates of~$U$,
the assertion easily transfers
from compact subsets of~$U$ to those of~$G$
(cf.\ the preceding proof).
\end{proof}
\section{\!\!Complex Kac-Moody groups and
their unitary forms}\label{SecKacMoody}
In this section, we generalize our results from Section~\ref{secphan}
about compact Lie groups | i.e.\ unitary forms of complex Lie groups
with respect to the compact involution | to unitary forms of
complex Kac-Moody groups with respect to the compact involution.
These groups correspond to complex
Kac-Moody algebras in very much the same way that
finite-dimensional semisimple complex Lie groups correspond to
finite-dimensional semisimple complex Lie algebras.
As far as complex Kac-Moody algebras are concerned,
we use the notations from~\cite{KacBook}:
If $A$ is a generalized Cartan matrix then
$\Lg(A)$ is the associated Kac-Moody Lie algebra and $\Lg'(A)$ is
its derived Lie algebra. The root space decomposition $\Lg(A) =
\Lh \oplus \bigoplus \Lg_\alpha$ induces a decomposition $\Lg'(A)
= \Lh' \oplus \bigoplus \Lg_\alpha$ and by \emph{loc.cit.}, Chapter~5,
the set of roots~$\tri$ decomposes into the set $\tri^{re}$ of real
roots and the set $\tri^{im}$ of imaginary roots.
If $W$ denotes the Weyl group of $\Lg(A)$ then the real roots are
canonically identified with the roots of a certain Coxeter system $(W,S)$
defined as in \emph{loc.cit.}, 3.13. For any real
root~$\alpha$, the root space $\Lg_\alpha$ is integrable.\\[3mm]
{\bf Convention:} We assume throughout that the Coxeter system
$(W,S)$ associated with $A$ is two-spherical and of finite rank
$N := |S|$.\\[3mm]
We now compile definitions, notation and facts
concerning complex Kac-Moody groups which will be used in the sequel.
For proofs, see~\cite{Remy} and the references therein.
By definition, the \emph{complex Kac-Moody group}
$G(A)$ is the group associated to the
(integrable representations of the) complex Lie algebra $\Lg'(A)$ in
the sense of \cite{KacConstructingGroups}, where it is also shown
that this definition is equivalent to the one given in~\cite{KP83}.
(For a more algebraic construction of Kac-Moody
groups which works over arbitrary fields, see~\cite{TitsFunctor}.)
In particular, if $F_{\Lg'(A)}$ denotes the set of
$\ad$-locally finite
elements in $\Lg'(A)$, then there exists a partial exponential
function $\exp\colon  F_{\Lg'(A)} \ra G(A)$, which is uniquely
determined by the property that for any integrable
$\Lg'(A)$-module $(V, d\pi)$, there exists
a representation $\pi$ of $G(A)$ on $V$
such that for any $v \in V$ the equality $(\pi \circ \exp(X)).v =
\sum_{n=0}^\infty \frac {(d\pi(X))^n.v}{n!}$ holds.
For any real root $\alpha \in \tri^{re}$,
we define the associated \emph{root group} to be
$U_\alpha := \exp{\Lg_\alpha}$. If $H := \bigcap_{\alpha\in\tri^{re}}
N_{G(A)}(U_\alpha)$, then the triple $(G, (U_\alpha)_{\alpha \in
\tri^{re}}, H)$ forms a root group data (RGD) system
(\emph{donn\'{e}e radicielle jumel\'{e}e}) of type $(W,S)$ in the sense
of \cite[1.5]{Remy}.
In particular, if $U_\pm := \langle
U_\alpha\,|\, \alpha \in \tri^{re}_\pm \rangle$,
then $B_\pm := HU_\pm$ are opposite Borel groups in $G(A)$
and one has the Bruhat and Birkhoff decompositions
\[
G(A) = \bigcup_{w \in W}B^+wB^+
=\bigcup_{w \in W} B^-wB^- = \bigcup_{w \in W} B^+wB^-.
\]
For a real root $\alpha$ the Lie algebra $\Lg(\alpha)$ generated
by $\Lg_\alpha$ and $\Lg_{-\alpha}$ in $\Lg'(A)$ is isomorphic to
$\mathfrak{sl}_2$ and contained in $F_{\Lg'(A)}$.
Therefore the inclusion map
$\phi_\alpha\colon \mathfrak{sl}_2 \cong  \Lg(\alpha)\hookrightarrow \Lg'(A)$
is integrable in the sense of
\cite{KacConstructingGroups} and thus gives rise to a map
$\phi_\alpha\colon \SL_2(\C) \ra G(A)$, whose image is denoted~$G_\alpha$.
It is known that $G_\alpha = \langle U_\alpha, U_{-\alpha} \rangle$
and that $\phi_\alpha$ is an isomorphism onto its image.
(See \cite{KP83} for details.)
We call $G_\alpha$ the \emph{rank one subgroup} associated
with the real root $\alpha$. 
If $\alpha_1, \dots, \alpha_n$ are real roots, then we define the
associated \emph{rank $n$ subgroup} to be $G_{\alpha_1, \dots,
\alpha_n} := \langle G_{\alpha_1}, \dots, G_{\alpha_n} \rangle$.
The rank one subgroups can be used
to define a topology
on $G(A)$ as was first observed by Kac and
Peterson in \cite{KP83}:
\begin{defn}
The final group topology with respect to the inclusion maps
$\phi_\alpha\colon \SL_2(\C) \ra G(A)$, where $\alpha$ runs through
the real roots of $\Lg(A)$ and $\SL_2(\C)$ is equipped with
its connected Lie group topology, is called the
\emph{Kac-Peterson topology}
on $G(A)$ and denoted $\tau_{KP}$.
\end{defn}
We record an easy fact for later reference:
\begin{la}\label{SimpleRoots}
If $\Pi = \{\alpha_1, \dots, \alpha_N\}$ is any choice of
simple roots for $\Lg(A)$, then $\tau_{KP}$ is the final group
topology with respect to the inclusions $\phi_{\alpha_i}\colon
\SL_2(\C) \ra G(A)$, $i=1, \dots, N$.
\end{la}
\begin{proof}
The Weyl group $W$ can be identified with $N_{G(A)}(H)/H$ and as
$H$ normalizes every root group, $W$ acts by conjugation on root groups.
This action is equivariant with respect
to the action on real roots, i.e.\ for any real root $\alpha$ and every
$w \in W$ one has $wU_\alpha w^{-1} = U_{w.\alpha}$ and thus
$wG_{\alpha}w^{-1} = G_{w.\alpha}$.
Now by definition every real root $\alpha$ of $\Lg'(A)$ is of the form
$w.\alpha_i$ for some $\alpha_i \in \Pi$ and thus every rank one subgroup
is conjugate to some $G_{\alpha_j}$.
Hence if $\tau$ is a group topology on $G(A)$ for which each of
the inclusions $\phi_{\alpha_i}\colon
\SL_2(\C) \ra G(A)$ is continuous for $i=1, \dots, N$,
then by continuity of conjugation all rank one subgroup inclusions
are in fact continuous. This implies the lemma.
\end{proof}
The following theorem relates complex Kac-Moody groups to
the main theme of this article:
\begin{thm}\label{KMKomega}
The topological group $(G(A), \tau_{KP})$ is a $k_\omega$-group.
\end{thm}
To prove Theorem \ref{KMKomega}, we need some preliminary
results on integrable representations.
\begin{la}\label{DifferentiableImpliesContinuous}
Let $(V,d\psi)$ be an integrable
$\g'(A)$-module
and $\psi\colon  G(A) \ra \Aut(V)$ be
the associated representation of~$G(A)$.
Let $\alpha$ be a real root.
Then every $v\in V$ is contained
in a finite-dimensional $G_\alpha$-submodule~$V_0$
of~$V$, and the orbit map
$G_\alpha \ra V$, $g \mapsto \psi(g).v$
is continuous into~$V_0$.
Moreover, for any $v\in V$ and real roots
$\alpha_1,\ldots, \alpha_n$, the map
\begin{equation}\label{betthere}
G_{\alpha_1}\times\cdots \times G_{\alpha_n}\to V\,,\quad
(g_1,\ldots, g_n)\mto \psi(g_1\cdots g_n).v
\end{equation}
has image in a finite-dimensional
vector subspace of~$V$
and is continuous into this space.
\end{la}
\begin{proof}
Since~$V$ is integrable and $\g_\alpha\sub F_{\g'(A)}$
is finite-dimensional,
$v$ is contained in a finite-dimensional
$\g_\alpha$-submodule~$V_0$ of~$V$
(see \cite[Proposition~3.6\,(a)]{KacBook} or part~(b) of the lemma
in \cite[\S1.2]{KacConstructingGroups}).
Then $\psi(\exp(X)).w=\sum_{n=0}^\infty\frac{(d\psi(X))^n.w}{n!}\in
V_0$ for all $X\in \g_\alpha$
and $w\in V_0$,
entailing that~$V_0$ is a $G_\alpha$-submodule
and that $\pi\colon G_\alpha\to \GL(V_0)$, $g\mto \psi(g)|_{V_0}$
is the smooth (and hence continuous)
representation of the Lie group~$G_\alpha$
with $d\pi(X)=d\psi(X)$ for $X\in \g_\alpha$.
As a consequence, the orbit map
$G_\alpha\to V_0$, $g\mto \psi(g).v=\pi(g).v$
is continuous.\\[2.5mm]
To prove the final assertion, we proceed
by induction. The case $n=1$
has just been settled.
As the inductive hypothesis,
we may now assume that
$\psi(G_{\alpha_2}G_{\alpha_3}\cdots G_{\alpha_n}).v$ is contained
in a finite-dimensional vector subspace
$V_0\sub V$ and that the map
$\phi\colon G_{\alpha_2}\times \cdots\times G_{\alpha_n}
\to V_0$,
$\phi(g_2,\ldots, g_n):=\psi(g_2\cdots g_n).v$
is continuous.
Each element of~$V$
being contained in a finite-dimensional
$\g_{\alpha_1}$-module,
the $\g_{\alpha_1}$-module
$V_1\sub V$ generated by~$V_0$
is finite-dimensional.
Now consider the map
$\pi\colon G_{\alpha_1}\to \GL(V_1)$,
$g\mto \psi(g)|_{V_1}$.
By the case $n=1$,
the map $G_{\alpha_1}\to V_1$, $g\mto \pi(g).w$
is continuous for each $w\in V_1$,
entailing that~$\pi$ is continuous.
The evaluation map
$\ve\colon \GL(V_1)\times V_1\to V_1$,
$(B,x)\mto B(x)$ being
continuous, it follows that
$\ve\circ (\pi\times \phi) \colon G_{\alpha_1}\times (G_{\alpha_2}\times
\cdots\times G_{\alpha_n})\to V_1$
is continuous. This is the mapping in~(\ref{betthere}).
\end{proof}
Now we can deduce Theorem \ref{KMKomega} as follows:\\[2.5mm]
{\bf Proof of Theorem \ref{KMKomega}.} Fix a system of
simple roots $\Pi = \{\alpha_1, \dots, \alpha_N\}$ and for every word
$I := (i_1, \dots, i_k)$ over $\{1, \dots, N\}$ denote by $G_I$ the
image of the map
\[
p_I\colon  G_{\alpha_{i_1}} \times \dots \times G_{\alpha_{i_k}} \ra G(A),
\quad (g_1, \dots, g_k) \mapsto g_1 \cdots g_k,
\]
which we equip with the quotient topology $\tau_I$ with respect to the
map $p_I$, where $G_{\alpha_j}$ is identified with $\SL_2(\C)$.
We claim that $G_I$ is a $k_\omega$-space for every~$I$.
As $\SL_2(\C)^k$ is a $k_\omega$-space, it suffices to show
(by Part (d) of Proposition \ref{basekomeg})
that $G_I$ is Hausdorff.
For this, let $g\not= h\in G_I$. 
By definition of $G(A)$, there
exists an integrable $\g'(A)$-module $(V,d\psi)$
such that $\psi(g)\not=\psi(h)$,
where $\psi\colon G(A)\to \Aut(V)$ is
the representation associated to~$d\psi$.
Thus there is $v\in V$ such that $\psi(g).v\not=\psi(h).v$.
Define $f \colon G_I\to V$,
$f(x):=\psi(x).v$.
The map $f \circ p_I$ is continuous
by Lemma~\ref{DifferentiableImpliesContinuous}.
Hence also~$f$ is continuous.
Now suitable preimages of~$f$
provide disjoint neighbourhoods of~$g$ and~$h$, finishing the proof
that $G_I$ is Hausdorff and $k_\omega$.
On the set $\mathcal W_N$ of words over $\{1, \dots, N\}$,
define a partial ordering by
putting $I \leq J$ if and only if
$I$ is a subsequence of $J$. Then for $I \leq J$ there is an
obvious inclusion map $G_I \ra G_J$.
Let $\tau_0$ denote the final topology with respect to the
system $(G_I, \tau_I)_{I \in \mathcal W_N}$ on the union
$G(A) = \bigcup G_I$. A
cofinal subsequence
$I_1 \leq I_2 \leq \cdots$ in $\mathcal W_N$ can be defined by
$I_1 := (1)$, $I_2 :=(1,2)$,
\dots, $I_{N+1} := (1,2, \dots, N,1)$ and so on.
We then have
\[
G_{I_1} \, \subseteq \, G_{I_2} \, \subseteq\,  \cdots \,
\subseteq \, G(A)
\]
and as each of the spaces $G_{I_n}$ is $k_\omega$,
we conclude from Proposition \ref{DLlkom} that
$(G(A), \tau_0)$ is a $k_\omega$-space.
It is easy to check that $\tau_0$ is a group
topology: The concatenation map
$G_{(i_1, \dots, i_k)} \times G_{(j_1, \dots, {j_m})}
\ra G_{(i_1, \dots, j_m)}$ is continuous and thus
$G_{(i_1, \dots, i_k)} \times G_{(j_1,
\dots, {j_m})} \ra G(A)$ is continuous.
Then by Proposition \ref{herecompat}, multiplication is continuous
and the proof for continuity of inversion is similar.
This shows that
$(G(A), \tau_0)$ is a $k_\omega$-group.
The inclusion maps $\phi_{\alpha_i}\colon
\SL_2(\C) \ra (G(A), \tau_0)$, $i=1, \dots, N$
are continuous by construction.
On the other hand, if $\tau_1$ is any group topology on $G(A)$
for which these inclusions are
continuous, then also the maps
$p_I\colon  G_{\alpha_{i_1}} \times \dots \times
G_{\alpha_{i_k}} \ra G(A)$
are continuous. In view of Lemma~\ref{SimpleRoots}, this implies the
theorem.\,\Punkt\\

\noindent We record a byproduct of the proof of Theorem~\ref{KMKomega}.
\begin{cor}\label{RFKomega}
Let $(V,d\psi)$ be an integrable
$\g'(A)$-module
and $\psi\colon  G(A) \ra \Aut(V)$ be
the associated representation of~$G(A)$.
Then $\psi$ is continuous with respect to the
Kac-Peterson topology on $G(A)$ and the topology of pointwise convergence
on $Aut(V)$.
\end{cor}
\begin{proof}
For each of the sets $G_I$, $I \in \mathcal W_N$ as defined in the
proof of Theorem~\ref{KMKomega}, the restriction $G_I \ra \Aut(V)$
of~$\psi$ is continuous by 
Lemma~{\rm\ref{DifferentiableImpliesContinuous}}.
As $G(A)=\dl\,G_I$,\vspace{-.8mm} the assertion follows.
\end{proof}
\noindent
Let us call a subalgebra $\Lh$ of $\Lg'(A)$ \emph{integrable}
if every ${\rm ad}_{\Lh}$-locally finite element of $\Lh$ is
${\rm ad}_{\Lg'(A)}$-locally finite.
(For a complex subalgebra $\Lh$ this is equivalent to the condition
that the inclusion map $\Lh \ra \Lg'(A)$ is integrable in the
sense of~\cite{KacConstructingGroups}.)
As usual, a subalgebra $\Lg^0(A) \subset \Lg'(A)$ is
called a \emph{real form} if
$\Lg'(A) = \Lg^0(A) \oplus i \Lg^0(A)\cong \Lg^0(A) \otimes \C$.
A real form $\Lg^0(A)$ induces an involutive Lie algebra automorphism
$\sigma$ on $\Lg'(A)$ by
\[
\sigma(X_0 + iX_1) := X_0 - iX_1 \quad(X_0, X_1 \in \Lg^0(A)).
\]
We observe:
\begin{prop}\label{RFIntegrable}
Let $\Lg^0(A) \subset \Lg'(A)$ be a real form and
$\sigma\colon \Lg'(A) \ra\Lg'(A)$ the associated involution.
Then the following hold:
\begin{itemize}
\item[\rm(a)]
$\Lg^0(A) \subseteq \Lg'(A)$ is an integrable subalgebra.
\item[\rm(b)]
$\sigma$ is integrable and the induced involution
$\sigma\colon G(A) \ra G(A)$ on the group level is
continuous.
\end{itemize}
\end{prop}
\begin{proof}
(a) Let us denote by ${\rm ad}$ and ${\rm ad}_0$ the adjoint representations
of $\Lg'(A)$ and $\Lg^0(A)$, respectively.
Let $X \in \Lg^0(A)$ be ${\rm ad}_0$-locally finite and let
$Y \in \Lg'(A)$.
We may write $Y = Y_0 + iY_1$ with $Y_0, Y_1 \in \Lg^0(A)$.
By assumption, there exist finite-dimensional
${\rm ad}_0(X)$-invariant subspaces $V_0$ and $V_1$ of
$\Lg^0(A)$ containing $Y_0$ and $Y_1$, respectively.
We thus have $Y \in V_0 \oplus i V_1$ and for
$Z = Z_0 + i Z_1 \in V_0 \oplus i V_1$, we obtain
\[
{\rm ad}(X)(Z) = [X, Z_0 + iZ_1]
= [X, Z_0] + i [X, Z_1]
= {\rm ad}_0(X)(Z_0) + i{\rm ad}_0(X)(Z_1) \in V_0 \oplus i V_1,
\]
showing that $V_0 \oplus i V_1$ is a finite-dimensional
${\rm ad}(X)$-invariant subspace containing~$Y$.
As $Y \in \Lg'(A)$ was arbitrary, we conclude
that $X$ is ${\rm ad}$-locally finite, which proves (a).

(b) Suppose that
$X$ is $\ad$-locally finite. For $Y \in \Lg'(A)$,
let $V$ be a finite-dimensional $\ad(X)$-invariant subspace
of $\Lg'(A)$ containing $\sigma(Y)$.
Put $V' := \sigma(V)$.
We claim that $V'$ is $\ad(\sigma(X))$-invariant.
Indeed, let $Z' \in V'$ with $Z' = \sigma(Z)$ for some $Z \in V$.
Then
\[
\ad(\sigma(X))(Z')
= [\sigma(X),\sigma(Z)]
= \sigma([X, Z])
= \sigma(\ad(X)(Z)) \in \sigma(V) = V'.
\]
Moreover, we have
$Y = \sigma(\sigma(Y)) \in V'$
and thus $V'$ is a finite-dimensional $\ad(\sigma(X))$-invariant
subspace containing~$Y$.
As~$Y$ was arbitrary, $\sigma(X)$ is
$\ad$-locally finite and thus $\sigma$ is integrable.
The induced involution $\sigma\colon G(A) \ra G(A)$ on $G(A)$
is determined by the fact that $\sigma(\exp(X)) = \exp(\sigma(X))$.
This formula shows at once that the restriction of $\sigma$ to any
rank one subgroup is smooth, which implies that~$\sigma$ is continuous.
\end{proof}
If $\Lg^0(A) \subseteq \Lg'(A)$ is a real form with associated
involution $\sigma$, then the group $G^0(A) := G(A)^\sigma$
is called the \emph{real form} of $G(A)$ associated with
$\Lg^0(A)$.
We equip $G^0(A)$ with the topology induced by~$G(A)$.
\begin{cor}\label{ReFoKomega}
Every real form of a complex Kac-Moody group is a $k_\omega$-group.
\end{cor}
\begin{proof}
Combine Proposition~\ref{RFIntegrable}\,(b)
with Proposition~\ref{basekomeg}\,(b)
and Theorem \ref{KMKomega}.
\end{proof}
In view of Part~(a) of Proposition~\ref{RFIntegrable},
one can think of another approach towards real forms of
complex Kac-Moody groups: The general construction
from~\cite{KacConstructingGroups} associates to
(the integrable complex representations of) any complex Lie
algebra~$\Lh$ an associated group~$H$.
The same construction can also be
applied to the integrable representations of any real Lie algebra.
In particular, we can use it
to associate a group~$G^0$ to any given real form
$\Lg^0(A) \subseteq \Lg'(A)$.
In view of Part~(a) of Proposition~\ref{RFIntegrable},
the inclusion $\Lg^0(A) \ra \Lg'(A)$ induces a map
$i\colon G^0 \ra G(A)$. Then one has:
\begin{prop}\label{RealForms}
The inclusion $\Lg^0(A) \ra \Lg'(A)$ induces an injection
$i\colon  G^0 \ra G^0(A) \subseteq G(A)$,
where $G^0(A)$ is the real form of $G(A)$
associated with $\Lg^0(A)$.
\end{prop}
\begin{proof}
Let $g \in G^0$. If $g \neq e$,
there exists an integrable
representation \mbox{$d\phi\colon \Lg^0(A) \ra \End(V^0)$}
with associated group representation
$\phi\colon G^0(A) \ra {\rm Aut}(V^0)$ such that
$\phi(g) \neq {\rm Id}_{V^0} = \phi(e)$.
We can complexify $d\phi$ to a Lie algebra homomorphism
$d\phi_\C\colon  \Lg^0(A) \otimes \C \ra \End(V^0 \otimes
\C)$ by $d\phi_\C(X \otimes z)(v) := d\phi(X)(zv)$,
which we may identify with a Lie algebra homomorphism
$d\phi_\C\colon  \Lg'(A) \ra \End(V^0 \otimes \C)$
by identifying $\Lg^0(A) \otimes \C$ with $\Lg'(A)$.
It is easy to see that $d\phi_\C$ is also integrable.
We thus obtain a commuting square
\[\begin{xy}
  \xymatrix{
      G(A) \ar^{\phi_\C}[r]  & \Aut(V^0 \otimes \C)\\
      G^0(A) \ar[u]^i  \ar^{\phi}[r]             &  \Aut(V^0) \ar[u]
  }
\end{xy},
\]
where $\Aut(V) \ra \Aut(V \otimes \C)$ is induced from the
inclusion $V \ra V \otimes \C$, $v \mapsto v \otimes 1$.
As $\phi(g) \neq {\rm Id}_{V^0}$, we have also
$\phi_{\C}(i(g)) \neq {\rm Id}_{V^0\otimes \C} = \phi_\C(i(e))$
and thus $i(g) \neq i(e)$, showing that $i$ is injective.
Since~$G^0$ is generated
by the elements
$\exp(X)$ for $X \in \Lg^0(A)$, which are $\sigma$-invariant,
we see that $i(G^0) \subseteq G^0(A)$.
\end{proof}
Proposition \ref{RealForms} suggests the question
of whether
the map $i$
actually defines
an isomorphism
$G^0 \cong G^0(A)$. This is true for the unitary real form
(defined below) and the split
real form (see \cite{Remy}). Moreover, in these cases $i(G^0) = G^0(A)$ is
generated by the groups $G^0_\alpha(A) := G^0(A) \cap G_\alpha$.
We do not know whether either of these statements remains
true for arbitrary real forms.\\[2.5mm]
In Section \ref{secphan} we obtained a presentation
of compact Lie groups, i.e.\ unitary real forms of complex semisimple
Lie groups, as Phan amalgams.
Recall that a unitary
form of a complex semi-simple Lie group~$G$ is the fixed point set
$K = G^\tau$ of the unique involution~$\tau$ on~$G$ whose
derivative
restricts to negative conjugated transpose on each rank~1
subalgebra.
Such an involution can also be defined for complex
Kac-Moody groups, see~\cite{Kac/Peterson:1985}.
We recall the construction here
to fix our notation:
Let~$\omega$ be the Chevalley involution of $\Lg(A)$
(see \cite[p.\,7]{KacBook}) and let $\overline \omega$ be the
\emph{compact involution} $\overline \omega (X) :=
\overline{\omega(X)}$.
This involution restricts to an involution
$\overline \omega$ of $\Lg'(A)$ and the subalgebra
$\Lk(A) := \{X \in \Lg'(A) \colon \overline{\omega}(X) = X\}
\subseteq \Lg'(A)$ is a real form (with associated involution
$\overline{\omega}$).
As $\overline \omega$ restricts to negative conjugated transpose
on each of the rank one subalgebras $\Lg(\alpha)$, we may call
$\Lk(A)$ the \emph{unitary form} of $\Lg'(A)$ in accordance
with the notion for complex semisimple Lie algebras.\\

\noindent The real form $K(A)$ of $G(A)$ corresponding to the unitary
form $\Lk(A)$ is also called the \emph{unitary form} of $G(A)$.
The involution
$\overline \omega\colon  G(A) \ra G(A)$
obtained from integrating $\overline \omega$ fixes every rank one
subgroup and satisfies $K_\alpha := G_\alpha^{\overline \omega}
= G_\alpha \cap K(A)$.
Moreover, $K(A)$ is generated by the groups
$K_\alpha$ for real $\alpha$. If we restrict the inclusion maps
$\phi_\alpha\colon  \SL_2(\C) \ra G_\alpha$ to
$\SU_2(\C)$, then we obtain isomorphisms
$\psi_\alpha\colon  \SU_2(\C) \ra K_\alpha$.
By Corollary~\ref{ReFoKomega}, the Kac-Peterson topology on $G(A)$
induces a $k_\omega$-topology on $K(A)$ which we also call
the \emph{Kac-Peterson topology on $K(A)$}
and denote by the same letter $\tau_{KP}$.
This topology can again be characterized by a universal property:
\begin{prop}\label{KPfinal}
If $\Pi = \{\alpha_1, \dots, \alpha_N\}$ is any choice of simple roots
for $\Lg(A)$, then the Kac-Peterson topology $\tau_{KP}$ on $K(A)$
is the final group topology with respect to the
inclusions $\psi_{\alpha_i}\colon
\SU_2(\C) \ra K(A)$, $i=1, \dots, N$.
\end{prop}
\begin{proof}
Equip the groups $K_{\alpha_i}$ with their compact connected
Lie group topology and let $\tau_0$ be the final group topology
with respect to the maps $\psi_{\alpha_i}$ or,
equivalently, the inclusions $K_{\alpha_i} \hookrightarrow K(A)$.
Note that the maps $\psi_{\alpha_i}$ are continuous
with respect to $\tau_{KP}$ and thus $\tau_0$ is finer than
$\tau_{KP}$;
in particular, $\tau_0$ is Hausdorff.
Similarly as in the proof of Theorem~\ref{KMKomega},
we associate to each word $I= (i_1, \dots, i_k)$ over
$\{1, \dots, N\}$ a map
\[
\psi_I\colon  K_{\alpha_{i_1}} \times \dots \times K_{\alpha_{i_k}}
\ra K(A), \quad (g_1, \dots, g_k) \mapsto g_1 \cdots g_k,
\]
whose image is denoted~$K_I$.
The maps $\psi_I$ are continuous both with respect to~$\tau_0$
and~$\tau_{KP}$.
As these topologies are Hausdorff and
$K_{\alpha_{i_1}} \times \dots \times K_{\alpha_{i_k}}$ is compact,
$\psi_I$ is a quotient map for both $\tau_0$ and $\tau_{KP}$.
In particular, $\tau_0$ and $\tau_{KP}$  coincide on each of the
subsets~$K_I$.
As in the proof of Theorem~\ref{KMKomega},
we see that $(K(A), \tau_0) = \dl\,K_I$.\vspace{-.6mm}
On the other hand, $(K(A), \tau_{KP})
= \dl\,(G_I \cap K(A))$.\vspace{-.6mm}
As $K_I \subseteq G_I \cap K(A)$
and the two topologies coincide on every~$K_I$, it
suffices to show that for each word $I \in \mathcal W_N$,
there exists $J \in \mathcal W_N$ such that $G_I \cap K(A) \subseteq K_J$.
This is essentially contained in the proof of Proposition~5.1
of~\cite{Kac/Peterson:1985}
but not stated explicitly, whence we elaborate
the argument here:
Denote by $s_j$ the root reflection at $\alpha_j$.
Then $S = \{s_1, \dots, s_N\}$
is a system of Coxeter generators for $W$.
Let us call a word $I=(i_1, \dots, i_k)$ \emph{special}
if $(s_{i_1}, \dots, s_{i_k})$ is
a reduced expression for the Weyl group element
$w(I) := s_{i_1} \cdots s_{i_k} \in W$.
Let us compute $G_I \cap K(A)$ in this case:
By the Bruhat decomposition for $\SL_2(\C)$, we have
$G_{\alpha_i} \subseteq B^+ \cup B^+ s_i B^+$.
It then follows from the general theory of Tits systems that
\[
G_I \, \subseteq \, (B^+ \cup B^+ s_{i_1}B^+)
\cdots  (B^+ \cup B^+ s_{i_k}B^+) \; \subseteq 
\bigcup_{w \leq w(I)} B^+wB^+ \, ,
\]
where the inequality $w \leq w(I)$ is understood with respect to the
Bruhat ordering of the Coxeter system $(W,S)$. Then
\[
G_I \cap K(A) \; \subseteq
\bigcup_{w \leq w(I)} K(A) \cap B^+wB^+\,.
\]
By Part~(c) of Proposition~5.1 in \cite{Kac/Peterson:1985},
we have $K(A) \cap B^+wB^+ \subseteq K_I \cdot T$ for each
$w \leq w(I)$, where $T$ is the compact torus defined by $T
:= H \cap K(A)$ and $H = \bigcap_{\alpha\in\tri^{re}}
N_{G(A)}(U_\alpha)$ as above.
This shows that
\[
G_I \cap K(A) \, \subseteq \, K_I\cdot T.
\]
Now $T \subseteq K_1 \cdots K_N$ and thus if we define $J$
as the concatenation of $I$ with $(1, \dots, N)$, then
$G_I \cap K(A) \subseteq K_J$.
If $I \in \mathcal W_N$ is arbitrary,
then this argument can be modified as follows:
Let $\mathcal W_I$ denote the set of special subwords of~$I$.
This is a finite subset of $\mathcal W_N$ and thus there exists a
special word $I' \in \mathcal W_N$ such that
\[
(\forall \, {J \in \mathcal W_I}) \;\, \; \; w(I') \geq w(J).
\]
Then as before the theory of Tits systems implies 
\[
G_I \; \subseteq \bigcup_{w \leq w({I'})}B^+wB^+
\]
and thus if we define $J$ as the concatenation of~$I'$ with
$(1, \dots, N)$, the same argument as before shows
$G_I \cap K(A) \subseteq K_J$, finishing the proof.
\end{proof}
\noindent Let us now see that the rank two subgroups are sufficient
for a characterization of the unitary form $K(A)$.
Adapting \cite[Definition~4.3]{Gra}, we declare:
\begin{defn}
Let $\Delta$ be a two-spherical Dynkin diagram
with set of labels $\Pi$, i.e., a Dynkin diagram whose induced
subdiagrams on at most two nodes are spherical.
Let $G$ be a topological group.
A \emph{topological
weak Phan system of type $\Delta$ over~$\C$ in~$G$}
is a family $(\overline{K}_j)_{j\in \Pi}$
of subgroups of~$G$
with the following properties:
\begin{description}
\item[twP1]
$G$ is generated by $\bigcup_{j \in \Pi} \overline{K}_j$;
\item[twP2]
$\overline{K}_j$ is isomorphic as a topological group
to a central quotient of a simply connected,
compact, semisimple Lie group of rank one,
for each $j\in \Pi$;
\item[twP3]
For all $i\not= j$ in~$\Pi$,
the subgroup
$\overline{K}_{ij}:=\langle \overline{K}_i\cup \overline{K}_j\rangle$
of~$G$ is isomorphic as a topological group
to a central quotient of the simply connected,
compact, semisimple Lie group of rank two belonging to the rank two
subdiagram of $\Delta$ on the nodes $i$ and $j$;
\item[twP4]
$(\overline{K}_i, \overline{K}_j)$ or $(\overline{K}_j, \overline{K}_i)$
is a standard pair in $\overline{K}_{ij}$
as defined in the paragraph before Theorem \ref{borovoi},
for all $i\not= j$ in~$\Pi$.
\end{description}
\end{defn}
\begin{defn}\label{phanamaldef}
Let $\Delta$ be a two-spherical Dynkin diagram with set of types~$\Pi$.
An amalgam
$\cA =(\overline{K}_j)_{j \in {\Pi \choose 1} \cup {\Pi \choose 2}}$
is called a \emph{Phan amalgam of type $\Delta$ over $\CC$},
if for all $\alpha, \beta \in \Pi$ the group
$\overline{K}_{\alpha\beta}$ is a central quotient of a simply connected
compact semisimple Lie group of type indicated by $\alpha$, $\beta$,
and if its morphisms are maps
$\overline{K}_\alpha \hookrightarrow \overline{K}_{\alpha\beta}$ which embed
$(\overline{K}_\alpha,\overline{K}_\beta)$
as a standard pair into $\overline{K}_{\alpha\beta}$
for $\alpha \neq \beta$.   
\end{defn}
A weak Phan system clearly gives rise to a Phan amalgam.
The converse is in general not true.
However, if there exists an enveloping group~$G$ of a Phan amalgam~$\cA$ into
which $\cA$ embeds, then the image of~$\cA$ in~$G$ is a weak Phan
system of~$G$.
Phan amalgams admitting such an enveloping group are called
\emph{strongly noncollapsing}.\\[3mm]
We have encountered examples of weak Phan systems in Section~\ref{secphan}.
More generally, let $A$ be a generalized Cartan matrix of
two-spherical and irreducible type $(W,S)$ and finite rank $|S|$ and
consider the unitary form $K(A)$ of the complex Kac-Moody group $G(A)$,
which we turn into a $k_\omega$-group
via the Kac-Peterson topology $\tau_{KP}$.
A topological weak Phan system for $(K(A), \tau_{KP})$
can be constructed as follows:
The type of the Phan system is the Dynkin diagram $\Delta$ of type $A$.
For any fundamental root $\alpha \in \Pi$ we put $\overline{K}_\alpha
:= K_{\alpha}$, which is identified as a topological group
with $\SU_2(\C)$ as above.
By construction, $(\overline{K}_\alpha)_{\alpha \in \Pi}$ is a topological
weak Phan system for $G$ and $(K(A),\tau_{KP})$ is an enveloping group of the
corresponding topological Phan amalgam.
This topological Phan amalgam is called the
\emph{topological standard Phan amalgam defined by $(K(A),\tau)$}.
We will show later that, if the Dynkin diagram $\Delta$ of $K(A)$ is a tree,
then there essentially exists a unique Phan amalgam associated to~$\Delta$.
On the other hand, if $\Delta$ admits cycles, this cannot be expected
in view of \cite[Section~6.5]{TitsFunctor} and~\cite{Muehlherr:1999}.\\[3mm]
We want to prove that any topological group $G$ with a topological
weak Phan system
of type~$\Delta$ over~$\C$, where~$\Delta$ is a tree,
carries a refined $k_\omega$-topology $\tau'$ such that
$(G,\tau')$ is a central
quotient of some $(K(A),\tau_{KP})$. As a first step in this direction,
we show that $(K(A),\tau_{KP})$ is a colimit of the amalgam defined by
its topological Phan system.
\begin{thm}\label{TopCompletion}
Let $K(A)$ the unitary form of some complex Kac-Moody group $G(A)$
associated with a symmetrizable generalized Cartan matrix of finite size
and two-spherical type $(W,S)$,
let $\tri^{re}$ be the set of real roots, and~$\Pi$ be a system
of fundamental roots.
Let $\I$ be the small category with objects
${\Pi\choose 1} \cup {\Pi \choose 2}$ and morphisms
$\{\alpha\} \to \{\alpha,\beta\}$,
for all $\alpha, \beta \in \Pi$, and let $\delta \colon \I \to \KOG$
be a diagram with $\delta(\{\alpha\}) = K_\alpha$,
$\delta(\{\alpha,\beta\}) = K_{\alpha\beta}$ and
$\delta(\{\alpha\} \to \{\alpha,\beta\})
= (K_\alpha \hookrightarrow K_{\alpha\beta})$.\\[2.5mm]
Then $((K(A), \tau_{KP}), (\iota_i)_{i \in { \Pi \choose 1}
\cup {\Pi \choose 2}})$, where
$\iota_i$ is the natural inclusion map, is a colimit of $\delta$
\begin{enumerate}
\item[{\rm (a)}] in the category $\G$ of abstract groups;
\item[{\rm (b)}] in the category $\HTG$ of Hausdorff topological groups;
\item[{\rm (c)}] in the category $\KOG$ of $k_\omega$ groups.
\end{enumerate}
\end{thm}
\begin{proof}
(a) is proved exactly as Theorem \ref{borovoi},
(b) follows from Proposition \ref{KPfinal} and Part~(c) follows from~(b)
and Corollary~\ref{colimkom}.
\end{proof}
We now state the main result of this section.
It is a general topological version of the main result of \cite{Gra}.
\begin{thm}[Topological Phan-type Theorem]\label{mainphan}
\,Let $\Delta$ be a two-spherical finite\linebreak
Dynkin tree-diagram
with generalized Coxeter matrix $A$
and let~$(G,\cO)$ be a Hausdorff topological
group admitting a topological
weak Phan system of type~$\Delta$ over~$\C$.
Then the topology on~$G$ can be refined to
a group topology making~$G$ a $k_\omega$-group isomorphic
to a central quotient of the unitary form $(K(A),\tau_{KP})$ of a complex
Kac-Moody group $G(A)$ equipped with its Kac-Peterson topology.

If~$\Delta$ is spherical, then the refined topology coincides
with the original topology, and $(G,\cO)$ is a compact
Lie group isomorphic to a central quotient of a simply connected\linebreak
compact semisimple Lie group whose universal complexification
is a simply connected\linebreak
complex semisimple Lie group of type~$\Delta$.
\end{thm}
Notice that we do \emph{not} claim that the group~$G$ in
Theorem~\ref{mainphan} is a topological quotient of the model $(K(A),\tau_{KP})$,
or a $k_\omega$-group. Both statements are in general
only true after refining the topology of $G$.
To illustrate the necessity of this refinement we consider
affine Kac-Moody groups: 
\begin{ex}
Let $\Delta_0$ be a spherical Dynkin diagram
and $A_0$ the associated Coxeter matrix. Then $G_0 := G(A_0)$ is a
finite-dimensional complex Lie group. Let us assume that~$G$ is simple
and not of type~$A_n$. Then the affine extension~$\Delta$
of~$\Delta_0$ is a tree and thus Theorem~\ref{mainphan} applies.
Denote by~$A$ the generalized Coxeter matrix of~$\Delta$ and
abbreviate $G := G(A)$, $K :=K(A)$.
We can consider~$G_0$ as an algebraic group over~$\C$ and define
the group $L_0G_0 := G_0(\C[[t]][t^{-1}])$ of $\C[[t]][t^{-1}]$-rational
points (cf. \cite[Section 13.2]{Kumar}).
The natural $\C^\times$-action on $\C[[t]][t^{-1}]$
by $(z.f)(t) := f(zt)$ induces a homomorphism\linebreak
$\C^\times \ra {\rm Aut}({L}_0G_0)$
and we write $\widetilde{L}_0G_0$ for the corresponding semidirect
product. By \cite[Definition 7.4.1 and Theorem 13.2.8]{Kumar},
there is a group homomorphism $q_0\colon  G \ra
\widetilde{L}_0(G_0)/C_0$ with central kernel,
where~$C_0$ denotes the center of $\widetilde{L}_0G_0$.
(Notice that our group~$G$ coincides with Kumar's
$\mathcal{G}_{\rm min}$).
Now $L_0G_0$ embeds naturally into the group
$LG_0 := C^\infty(\mathbb \C^\times, G_0)$
and similarly $\widetilde{L}_0G_0$ embeds into
\[
\widetilde{L}G_0 :=  C^\infty(\mathbb \C^\times, G_0) \ltimes \C^\times,
\]
where the $\C^\times$-action is given by precomposition.
If we equip $C^\infty(\mathbb \C^\times, G_0)$ with the topology
of uniform convergence of all derivatives and $\widetilde{L}G_0$
with the corresponding product topology, then the latter becomes a
metrizable topological group. The same is true of the quotient
$\widetilde{L}_0(G_0)/C$ by the center~$C$.
Denote by $(\widehat{K}, \tau)$ the image of~$K$ under the
homomorphism
\[
q\colon  G \xra{q_0} \widetilde{L}_0(G_0)/C_0 \hookrightarrow \widetilde{L}G_0/C
\]
equipped with the topology induced from the metrizable topology
on $\widetilde{L}G_0/C$ introduced above.
Denote by $(K_j)_{j \in \Pi}$ the standard Phan system in~$K$
and write $\overline{K_j}$ for the image of~$K_j$ under~$q$.
We observe that~$q|_{K_j}$ is continuous and that $\overline{K_j}$
is a central quotient of~$K_j$ (since~$q$ has central kernel).
These facts combine to show that $(\overline{K_j})_{j \in \Pi}$
is a topological weak Phan system for $(\widehat{K}, \tau)$.
On the other hand, since $\widehat{K}$ is metrizable and
not locally compact, it cannot be~$k_\omega$ by
Proposition~\ref{Metrizable}.
However, the topology~$\tau$ is easily refined to
a $k_\omega$-topology~$\tau_0$ (we can take the
topology making~$q|_K$ a quotient map).
\end{ex}
Notice that one can produce a lot of non-$k_\omega$ topologies
with weak topological Phan systems on the group $\widehat K$ by using
various topologies on $C^\infty(\mathbb \C^\times, G_0)$.
These topologies also lead to a multitude of topologies on the flag
variety (or building) of $G(A)$ on which $\widehat K$ acts
transitively.
We would like to advertise the topology induced from the
$k_\omega$-topology on $\widehat K$ (or, equivalently, $G(A)$)
as the most canonical topology on the flag variety.
It has the\linebreak
important advantage that the Schubert decomposition of
the flag variety becomes an honest CW decomposition.
Moreover, with this topology the geometric realization of the
associated topological building will be contractible.
The latter property is a crucial ingredient in the Kramer-Mitchell-Quillen
proof of Bott periodicity (\cite{Kramer}, \cite{Mitchell}).
Finally, let us mention that the $k_\omega$-topology
that we call Kac-Petersen topology occurs under various names
in the literature. In~\cite{Kumar}, it is called the
\emph{analytic topology} (more precisely, the analytic topology
on Kumar's group $\mathcal G$ restricts to the Kac-Petersen topology
on our group $G(A) \subseteq \mathcal G$), to distinguish it
from the Zariski topology.
\begin{defn} \label{3.4}
Let
$\mathcal{A} = (P_1 \stackrel{\iota_1}{\hookleftarrow} (P_1 \cap P_2)
\stackrel{\iota_2}{\hookrightarrow} P_2)$
and
$\mathcal{A}' = (P'_1 \stackrel{\iota_1'}{\hookleftarrow}
(P'_1 \cap P'_2) \stackrel{\iota_2'}{\hookrightarrow} P'_2)$
be amalgams consisting of abstract groups.
The amalgams $\mathcal{A}$ and $\mathcal{A}'$ are of the same
\emph{type} if there exist isomorphisms
$\phi_i \colon P_i \rightarrow P_i'$ such that
$(\phi_i \circ \iota_i)(P_1 \cap P_2) = \iota_i'(P_1' \cap P_2')$
for $i = 1,2$.
If, in addition, there exists an isomorphism
$\phi_{12} : P_1 \cap P_2 \to P_1' \cap P_2'$ such that the diagram
\[
\xymatrix{
 & P_1 \ar[rrr]^{\phi_1} & & & P_1'\\
P_1 \cap P_2 \ar[ur]^{\iota_1} \ar[dr]_{\iota_2} \ar[rrr]^{\phi_{12}}
& & & P_1' \cap P_2' \ar[ur]^{\iota_1'} \ar[dr]_{\iota_2'} & \\
 & P_2 \ar[rrr]^{\phi_2} & & & P_2'
}
\]
commutes,
then the amalgams $\mathcal{A}$ and $\mathcal{A}'$ are called
\emph{isomorphic}. For sake of brevity, we denote an isomorphism between
the amalgams $\mathcal{A}$ and $\mathcal{A}'$ by the triple
$(\phi_1, \phi_{12}, \phi_2)$. \\ [2mm]
In case the amalgams consist of topological groups, the
isomorphisms involved are assumed to be topological.
\end{defn}
\begin{la}[Goldschmidt's Lemma; cf.\
\cite{Goldschmidt:1980}] \label{Goldschmidt}
Let $\mathcal{A} = (P_1 \hookleftarrow (P_1 \cap P_2) \hookrightarrow P_2)$
be an amalgam consisting of topological groups,
let $A_i = \mathrm{Stab}_{\Aut(P_i)}(P_1 \cap P_2)$ for $i = 1,2$,
and let $\alpha_i \colon  A_i \rightarrow \Aut(P_1 \cap P_2)$
be homomorphisms mapping $a \in A_i$ onto its restriction to
$P_1 \cap P_2$.
Then there is a one-to-one correspondence between isomorphism classes
of amalgams of the same type as $\mathcal{A}$ and
$\alpha_2(A_2)$-$\alpha_1(A_1)$ double cosets in $\Aut(P_1 \cap P_2)$.
In other words, there is a one-to-one correspondence between the
different isomorphism types of amalgams
$P_1 \hookleftarrow (P_1 \cap P_2) \hookrightarrow P_2$ and the
double cosets $\alpha_2(A_2) \backslash \Aut(P_1 \cap P_2) / \alpha_1(A_1)$.
\end{la}
\begin{proof}
By Definition \ref{3.4}, any amalgam of the same type as
$P_1 \stackrel{\iota_1}{\hookleftarrow}
(P_1 \cap P_2) \stackrel{\iota_2}{\hookrightarrow} P_2$
is isomorphic to an amalgam of the form
$P_1 \stackrel{\iota_1 \circ \beta_1}{\hookleftarrow}
(P_1 \cap P_2) \stackrel{\iota_2 \circ \beta_2}{\hookrightarrow} P_2$,
where $\beta_1$ and $\beta_2$ are automorphisms of $P_1 \cap P_2$.
Such an amalgam is isomorphic to
$P_1 \stackrel{\iota_1}{\hookleftarrow} (P_1 \cap P_2)
\stackrel{\iota_2 \circ \beta_2 \circ \beta_1^{-1}}{\hookrightarrow} P_2$
by the amalgam isomorphism $(\mathrm{id}, \beta_1, \mathrm{id})$.
It remains to decide when two amalgams
$P_1 \stackrel{\iota_1}{\hookleftarrow}
(P_1 \cap P_2) \stackrel{\iota_2 \circ \gamma_1}{\hookrightarrow} P_2$
and $P_1 \stackrel{\iota_1}{\hookleftarrow}
(P_1 \cap P_2) \stackrel{\iota_2 \circ \gamma_2}{\hookrightarrow} P_2$
are isomorphic. That is, we want to find an amalgam isomorphism
$(\phi_1, \phi_{12}, \phi_2)$ such that
$\iota_1 \circ \phi_{12} = \phi_1 \circ \iota_1$ and
$\iota_2 \circ \gamma_2 \circ \phi_{12} = \phi_2 \circ \iota_2
\circ \gamma_1$.
The former equality says
$\phi_{12} \in \alpha_1(A_1)$ while the latter equality says
$\gamma_2 \circ \phi_{12} \circ \gamma_1^{-1} \in \alpha_2(A_2)$.
Therefore an isomorphism between these two amalgams exists if and only if
$\alpha_1(A_1) \cap \gamma_2^{-1} \alpha_2(A_2) \gamma_1 \neq \emptyset$.
The latter is equivalent to $\alpha_2(A_2) \gamma_1 \alpha_1(A_1)
= \alpha_2(A_2) \gamma_2 \alpha_1(A_1)$.
\end{proof}
\begin{defn}
A Phan amalgam $(\overline{K}_{\alpha\beta})_{\alpha, \beta \in \Pi}$,
and any subamalgam of a Phan amalgam, is called \emph{unambiguous}
if every $\overline{K}_{\alpha\beta}$ is isomorphic to the
corresponding $K_{\alpha\beta}$.
\end{defn}
\begin{prop} \label{refinement} \label{similar}
Every topological Phan amalgam
$\overline{\cA}=
(\overline{K}_J)_{J \in {\Pi \choose 1} \cup {\Pi \choose 2}}$
has an unambiguous covering
$\cA=({G}_J)_{J \in {\Pi \choose 1} \cup {\Pi \choose 2}}$
that is unique up to equivalence of coverings. 
Furthermore, every strongly noncollapsing Phan amalgam
$\overline{\cA}$ has a unique $($up to equivalence of coverings$)$
unambiguous strongly noncollapsing covering~$\cA$.
\end{prop}
\begin{proof} 
We proceed by induction on $|S|$, where $S$ is a subset of
${\Pi \choose 1} \cup {\Pi \choose 2}$ with the property that for each
$J \in S$, the inclusion $\emptyset \neq J' \subseteq J$ implies $J' \in S$.
We order $S$ by inclusion, i.e., the maximal elements of $S$ are those
sets which are not properly contained in another element of~$S$.
In particular, if~$J$ is a maximal element of~$S$, then
$S \backslash \{ J \}$ is again a subset of
${\Pi \choose 1} \cup {\Pi \choose 2}$ with the property that for
each $J' \in S \backslash \{ J \}$, the inclusion
$\emptyset \neq J'' \subseteq J'$ implies $J'' \in S \backslash \{ J \}$.
Let $\overline{\cA}_S=(\overline{K}_J)_{J \in S}$
with morphisms
$\overline{\iota}_{J',J} : \overline{K}_{J'} \to \overline{K}_{J}$
for $J' \subsetneq J$.
The case $S=\emptyset$, which vacuously yields an unambiguous amalgam,
serves as basis of the induction. \\[2.5mm]
Suppose that $S$ is non-empty and that for  every subset
$S'\subsetneq S$ the conclusion of the proposition holds.
Let $J$ be a maximal element of $S$,
define $S'=S\setminus\{J\}$, and consider
$\overline{\cA}_{S'}=(\overline{K}_{J'})_{J'\in S'}$.
By the inductive hypothesis,
there is a unique unambiguous covering 
amalgam $(\cA_{S'}=(K_{J'})_{J'\in S'},\pi')$ of $\overline{\cA}_{S'}$. 
By the definition of a Phan amalgam, there is a quotient map $\psi$ from 
$K_J$ onto $\overline{K}_J$. In case $|J| = 1$ define $\cA_S$ as the union
of $\cA_{S'}$ and $\{ K_J \}$ without additional morphisms.
Clearly, $\cA_S$ is an unambiguous covering of $\overline{\cA}_{S}$.
In case $|J| = 2$ assume $J = \{ \alpha, \beta \}$.
For $J' \subsetneq J$ let
$\iota_{J',J} : K_{J'} \hookrightarrow K_J$ be the natural inclusion map.
The function $\psi$ maps $\mc{G} = (\iota_{J',J}(K_{J'}))_{J'\subsetneq J}$
onto $(\overline{\iota}_{J',J}(\overline{K}_{J'}))_{J'\subsetneq J}$.
Indeed, the standard pair $(K_{\{ \alpha \}},K_{\{ \beta \}})$
of $K_{\{ \alpha, \beta \}}$ is mapped by~$\psi$ onto the standard
pair $(\overline{K}_{\{ \alpha \}},\overline{K}_{\{ \beta \}})$
of $\overline{K}_{\{ \alpha, \beta \}}$, cf.\ \ref{phanamaldef}
and the paragraph before~\ref{borovoi}.
Furthermore, by the
inductive hypothesis there is an amalgam isomorphism
$\phi \colon (\mc{G},\pi'|_{{\mc{G}}}) \to (\mc{G},\psi|_{\mc{G}})$
such that $\psi \circ \phi =\pi'|_{{\mc{G}}}$.
Therefore we can define an unambiguous covering~$\cA_S$
of $\overline{\cA}_{S}$ as the union of $\cA_{S'}$ and~$\{ K_J \}$
with the additional morphisms $\iota_{J',J} \circ \phi$ for $J' \subsetneq J$.
This completes the proof
of the existence of an unambiguous covering~$\cA_S$.\\[2.5mm]
It remains to prove the uniqueness of the unambiguous covering~$\cA_S$
that we have just constructed.
Let $\cB_S=(B_J)_{J\in S}$ and $\cC_S=(C_J)_{J\in S}$
be two such coverings with amalgam homomorphism~$\pi_1$,
respectively~$\pi_2$ onto~$\cA$.
By the inductive hypothesis, there exists an isomorphism
$\phi \colon \cB_{S'} \to \cC_{S'}$ with
$\pi_1|_{\cB_{S'}}=\pi_2 \circ \phi$.  
In order to extend $\phi$ to $B_J$, we have to deal with two cases:
First, let $J = \lb \alpha, \beta \rb$ where $\alpha$ and $\beta$
are orthogonal roots.  In this 
case, $B_{\alpha\beta} \cong C_{\alpha\beta} \cong K_{\alpha\beta}$
are isomorphic to a direct product of
$B_{\alpha} \cong C_{\alpha} \cong K_{\alpha}$ and 
$B_{\beta} \cong C_{\beta} \cong K_{\beta}$.  Clearly, 
$\phi$ is already known on $B_{\alpha}$ and $B_{\beta}$, and so $\phi$ 
extends uniquely to~$B_{\alpha\beta}$.
This extension, also denoted~$\phi$, is a 
well-defined amalgam isomorphism from~$\cB_S$ to~$\cC_S$, and furthermore
$\pi_1=\pi_2 \circ \phi$ holds.
In the second case, $B_J\cong C_J\cong K_J$ is isomorphic to a simply
connected compact almost simple Lie group of rank one or two.
By the universality of the covering
$\pi_1 \colon B_J \rightarrow \overline{K}_J$, as~$B_J$ is simply connected,
there exists a unique homeomorphism $\psi \colon B_J \rightarrow C_J$
such that $\pi_1 = \pi_2 \circ \psi$.
Consider the map $\alpha \colon \overline{K}_J \to \overline{K}_J$,
$u \mapsto (\pi_2 \circ \psi \circ \pi_1^{-1})(u)$. 
This map~$\alpha$ is a well-defined automorphism of~$\overline{K}_J$,
because the cosets of the kernel of~$\pi_1$ are mapped by~$\psi$ to cosets of 
the kernel of $\pi_2$.
Every  automorphism of $\overline{K}_J$, in particular~$\alpha$,
is continuous by a corollary of Got\^o's Commutator Theorem
(see \cite[Corollary 6.56]{Hofmann/Morris:1998})
and van der Waerden's Continuity Theorem (cf.\
\cite[Theorem 5.64]{Hofmann/Morris:1998}).
By~\cite{Kawada:1940}, $\alpha$
lifts to a unique continuous automorphism of~$C_J$
(see also~\cite{Hofmann:1962}).
In other words, there is a unique (continuous) automorphism~$\beta$
of~$C_J$ such that $\pi_2 \circ \beta=\alpha \circ \pi_2$.
Define $\theta \colon B_J \rightarrow C_J$,
$\theta(b) = (\beta^{-1} \circ \psi)(b)$.
First of all, by definition we have 
$\pi_1|_{B_J}=\pi_2 \circ \theta$,~as 
\begin{eqnarray*}
\pi_2 \circ \theta & = & \pi_2 \circ \beta^{-1} \circ \psi \\
& = & \alpha^{-1} \circ \pi_2 \circ \psi \\
& = & \pi_1|_{B_J} \circ \psi^{-1} \circ \pi_2^{-1}|_{\overline{K}_J}
\circ \pi_2 \circ \psi \\
& = & \pi_1|_{B_J}.
\end{eqnarray*}
Second, for every $J'\sub J$, we have 
that $\theta^{-1} \circ \phi_{|B_{J'}}$ is a lifting to $B_{J'}$ of the 
identity automorphism of $\overline{K}_{J'}$ and, by the above,
it is the identity.
This holds due to
$\theta^{-1} \circ \phi_{|B_{J'}}
= \psi^{-1} \circ \beta \circ \phi_{|B_{J'}}$
and the following, considered on $B_{J'}/\ker({\pi_1}_{|B_{J'}})$:
\begin{eqnarray*}
\psi^{-1}\circ{\pi_2}^{-1}_{|C_{J'}}\circ\alpha\circ\pi_2\circ\phi_{|B_{J'}}
& = &
\psi^{-1} \circ {\pi_2}^{-1}_{|C_{J'}} \circ \pi_2 \circ \psi \circ
{\pi_1}^{-1}_{|B_{J'}} \circ \pi_2 \circ \phi_{|B_{J'}} \\
& = & {\pi_1}^{-1}_{|B_{J'}} \circ \pi_2 \circ \phi_{|B_{J'}} \\
& = & \id . 
\end{eqnarray*}
Thus~$\phi$ and~$\theta$ agree on every subgroup $B_{J'}$, which 
allows us to extend~$\phi$ to all of~$\cB_S$ by defining it on~$B_J$
as~$\theta$.
Finally, if $\overline{\cA}$ is strongly noncollapsing, 
so is its unique unambiguous covering $\cA$, finishing the proof. 
\end{proof}
\begin{prop} \label{uniquecover}
Let $n \geq 2$, and let 
$\cA$ be a tree-like strongly noncollapsing unambiguous irreducible
Phan amalgam of rank~$n$.
Then~$\mc{A}$ is unique up to isomorphism, i.e.,
$\mc{A}$ is isomorphic to a standard Phan amalgam.
\end{prop}
\begin{proof}
Let $\cA=(K_J)_{J \in {\Pi \choose 1} \cup {\Pi \choose 2}}$
be a tree-like unambiguous strongly noncollapsing irreducible Phan 
amalgam of rank at least two.
We proceed by induction on~$n$.
The amalgams of rank two are unique by definition. 
\subsubsection*{Rank three}
Let $\Pi = \{ \alpha, \beta, \gamma \}$ and let~$\cA$ and~$\cA'$
be the amalgams 
\[
\xymatrix{
K_\alpha \ar[r] \ar[dr] & K_{\alpha\beta} \\
K_\beta \ar[ur] \ar[dr] & K_{\alpha\gamma} \\
K_\gamma \ar[ur] \ar[r] & K_{\beta\gamma}
} \quad \quad \mbox{respectively} \quad \quad
\xymatrix{
K'_\alpha \ar[r] \ar[dr] & K'_{\alpha\beta} \\
K'_\beta \ar[ur] \ar[dr] & K'_{\alpha\gamma} \\
K'_\gamma \ar[ur] \ar[r] & K'_{\beta\gamma}
}
\]
According to Goldschmidt's Lemma (Lemma \ref{Goldschmidt}),
the amalgams $\mc{B}=
(K_{\alpha\beta} \hookleftarrow K_{\beta} \hookrightarrow K_{\beta\gamma})$
and $\mc{B}'=
(K'_{\alpha\beta} \hookleftarrow K'_{\beta} \hookrightarrow K'_{\beta\gamma})$
are isomorphic via some amalgam isomorphism~$\psi$, because every
automorphism of the group~$K_{\beta}$ is induced by some automorphism
of the group $K_{\alpha\beta}$
(cf.\ \cite[Theorem 6.73]{Hofmann/Morris:1998}). 
So clearly $\psi(K_{\beta}) = K'_{\beta}$.\\[2.5mm]
The groups $K_{\alpha}$ and $K_{\beta}$ form a standard pair
in~$K_{\alpha\beta}$, and hence $\psi(K_{\alpha})$ and
$K'_{\beta}=\psi(K_{\beta})$ form a standard pair in
$K'_{\alpha\beta}=\psi(K_{\alpha\beta})$.
Therefore, since standard pairs are conjugate, there exists an automorphism
of $K'_{\alpha\beta}$ that maps $\psi(K_{\alpha})$ onto~$K'_{\alpha}$
and that leaves $K'_{\beta}=\psi(K_{\beta})$ invariant.
Thus, we can assume $\psi(K_{\alpha})=K'_{\alpha}$.\\[2.5mm]
Define $D_\alpha = N_{K_{\alpha}}(K_{\beta})$ and
$D_\gamma = N_{K_{\gamma}}(K_{\beta})$, where the groups~$K_{\beta}$,
$K_{\alpha}$ are considered as subgroups of~$K_{\alpha\beta}$ and
the groups $K_{\beta}$, $K_{\gamma}$ are considered as subgroups
of~$K_{\beta\gamma}$.
Since~$K_{\beta}$ and~$K_{\alpha}$ form a standard pair
in~$K_{\alpha\beta}$, the group~$D_\alpha$ is a maximal torus in~$K_{\alpha}$.
For the same reason, $D_\gamma$ is a maximal torus in~$K_{\gamma}$.
Moreover, let $D^\alpha_\beta = N_{K_{\beta}}(K_{\alpha})$ and
$D^\gamma_\beta = N_{K_{\beta}}(K_{\gamma})$.
Again, these are maximal tori in $K_{\beta}$,
satisfying $D^i_\beta = C_{K_{\beta}}(D_i)$ for $i = \alpha, \gamma$.
Let $\pi \colon \cA \to G$ be an enveloping group of~$\mc{A}$ such that~$\pi$
is injective on every $\overline{K}_{i}$, $i \in \Pi$, which is possible
since $\cA$ is strongly noncollapsing.
Then $\pi(D_\alpha^i) = C_{\pi(K_{\beta})}(\pi(D_i))$.
Moreover, since $\pi(K_{\beta})$ and $\pi(D_\gamma)$ are invariant
under $\pi(D_\alpha)=N_{\pi(K_{\alpha})}(\pi(K_{\beta}))$,
the group $\pi(D^\gamma_\beta)=C_{\pi(K_{\beta})}(\pi(D_\gamma))$
is invariant under $\pi(D_\alpha)$.
So the maximal torus $\pi(D_\alpha)$ of $\pi(K_{\alpha})$ leaves
invariant the maximal tori $\pi(D^\alpha_\beta)$ and
$\pi(D^\gamma_\beta)$ of $\pi(K_{\beta})$.
Analysis of the rank two group $\pi(K_{\alpha\beta})$ shows that
$D_\beta^\alpha = D^\gamma_\beta$. 
Hence we can use the notation $D_\beta=D_\beta^\alpha=D_\beta^\gamma$.
{}From $N_{K_{\beta}}(K_{\alpha}) = D_\beta$ it follows that
$\psi(D_\beta)=D_\beta' := N_{K'_{\beta}}(K'_{\alpha})
= N_{K'_{\beta}}(K'_{\gamma})$,
because $\psi(K_{\alpha})=K'_{\alpha}$ and $\psi(K_{\beta})=K'_{\beta}$.
Therefore we have identified a root system of type $A_2$, $B_2$,
or $G_2$ with respect to which the groups $K_\beta' = \psi(K_\beta)$,
$K'_\gamma$, $\psi(K_\gamma)$ occur as fundamental groups.
By inspection (cf.\ \cite[Plate X]{Bourbaki:2002}),
there exist automorphisms of $K'_{\beta\gamma}$ mapping $\psi(K_\gamma)$
onto $K'_\gamma$ and centralizing $K'_\beta$.
Hence we can assume $\psi(K_\gamma) = K'_\gamma$.\\[1.9mm]
Finally, we have $\psi(K_{\alpha\gamma}) = K'_{\alpha\gamma}$,
because $K_{\alpha\gamma} \cong K_\alpha \times K_\gamma$
and $K'_{\alpha\gamma} \cong K'_\alpha \times K'_\gamma$.\vspace{-1.2mm}
\subsubsection*{Rank at least four}
Let $|\Pi| \geq 4$ and let
$\cA=(K_J)_{J \in {\Pi \choose 1} \cup {\Pi \choose 2}}$
be a tree-like unambiguous strongly noncollapsing irreducible Phan 
amalgam.
Then there exists a unique amalgam $\cB_\cA$ obtained by adding a
colimit $H_\alpha$ of
$\cA_\alpha :=
(K_J)_{J \in {\Pi \backslash \{ \alpha \} \choose 1}
\cup {\Pi \backslash \{ \alpha \} \choose 2}}$
and a colimit $H_\beta$ of $\cA_\beta :=
(K_J)_{J \in {\Pi \backslash \{ \beta \} \choose 1}
\cup {\Pi \backslash \{ \beta \} \choose 2}}$
to the amalgam $\cA$ where $\alpha$ and $\beta$ are perpendicular in $\Pi$.
By the inductive hypothesis,
both $\cA_\alpha$ and $\cA_\beta$ are isomorphic to some standard Phan amalgam,
so the colimits $H_\alpha$ and $H_\beta$ are determined in
Theorem \ref{TopCompletion}.\\[2.5mm]
In particular, if~$\cA$ and~$\cA'$ are tree-like unambiguous
strongly noncollapsing irreducible Phan amalgams,
then $H_\alpha \cong H'_\alpha$ and $H_\beta \cong H'_\beta$.
For $\cB := (K_J)_{J \in {\Pi \backslash \{ \alpha, \beta \} \choose 1}
\cup {\Pi \backslash \{ \alpha, \beta \} \choose 2}}$
and $\cB' := (K'_J)_{J \in {\Pi \backslash \{ \alpha, \beta \} \choose 1}\cup
{\Pi \backslash \{ \alpha, \beta \} \choose 2}}$,
again by induction $\cB \cong \cB'$, so | for a colimit $H_0$ of $\cB$
and a colimit $H'_0$ of $\cB'$ | the amalgams
$H_\alpha \hookleftarrow H_0 \hookrightarrow H_\beta$
and $H'_\alpha \hookleftarrow H'_0 \hookrightarrow H'_\beta$ are of
the same type; the isomorphism type of $H_0 \cong H'_0$ again is
determined in Theorem~\ref{TopCompletion}.
By \cite[Theorem 8.2]{Caprace:2005} and Goldschmidt's Lemma
(Lemma \ref{Goldschmidt}), the amalgams
$H_\alpha \hookleftarrow H_0 \hookrightarrow H_\beta$
and $H'_\alpha \hookleftarrow H'_0 \hookrightarrow H'_\beta$
are isomorphic under some map $\phi$.
The standard Phan amalgams $\phi(\cB)$ and $\cB'$ in $H'_0$ correspond to
two choices of a maximal torus of $H'_0$, which are conjugate
by the Iwasawa decomposition.
So, correcting $\phi$ by an inner automorphism of~$H'_0$,
we may assume that $\phi(K_{i}) = K'_{i}$ for all
$i \in {\Pi \backslash \{ \alpha, \beta \} \choose 1} \cup
{\Pi \backslash \{ \alpha, \beta \} \choose 2}$.
Also, by studying  $H_\alpha'$ and $H_\beta'$, we have
$\phi(K_{\alpha}) = K'_{\alpha}$, $\phi(K_{\beta}) = K'_{\beta}$
and $\phi(K_{\alpha i}) = K'_{\alpha i}$,
$\phi(K_{\beta i}) = K'_{\beta i}$ for all
$i \in \Pi \backslash \{ \alpha, \beta \}$.
It remains to realize that $K_{\alpha\beta}$ is the direct product
of~$K_{\alpha}$ and~$K_{\beta}$, so that~$\phi$ induces an isomorphism
between~$\mc{A}$ and~$\mc{A}'$.\vspace{-.5mm}
\end{proof}
{\bf Proof of Theorem \ref{mainphan}.}
The weak Phan system of~$G$ gives rise to a strongly noncollapsing
Phan amalgam~$\mc{A}$, which by Proposition~\ref{similar}
is covered by a unique strongly noncollapsing unambiguous
Phan amalgam~$\widehat{\mc{A}}$.
This strongly noncollapsing unambiguous Phan amalgam~$\widehat{\mc{A}}$
is isomorphic to a standard Phan amalgam by Proposition~\ref{uniquecover}.
The claim then follows from Theorem~\ref{TopCompletion}
(and Theorem~\ref{borovoi}).\,\Punkt \\ [3mm]
Let us now return
to the Curtis-Tits Theorem
which motivated our research on topological presentations
resulting in this section of the
current article.
As we
mentioned in the introduction, the Curtis-Tits Theorem
was
originally proved
in \cite{Curtis:1965}, \cite{Tits:1962}.
By a result of Abramenko and M\"uhlherr \cite{Abramenko/Muehlherr:1997}
(or, alternatively, M\"uhlherr \cite{Muehlherr};
see also \cite{Bennett/Gramlich/Hoffman/Shpectorov:2003}, \cite{Gramlich}),
the Curtis-Tits Theorem follows from the simple connectedness
of the opposites geometry
that is associated to a twin building as follows.
Given a twin building $\mathcal{T}=(\mathcal{B}_+,\mathcal{B}_-,\delta_*)$
consisting of the buildings $\mathcal{B}_+ = (\mathcal{C}_+,\delta^+)$
and $\mathcal{B}_- = (\mathcal{C}_-,\delta^-)$ and the codistance
$\delta_*$ (cf.\ \cite{Tits:1992}), define the opposites
chamber system as $\mathrm{Opp}(\mathcal{T})
=\{(c_+,c_-)\in \mathcal{C}_+\times \mathcal{C}_- \mid\,
\delta_*(c_+,c_-)=1 \}$.
(Chambers $x\in\mathcal{C}_+$ and $y\in\mathcal{C}_-$ with
$\delta_*(x,y)=1$ are called \emph{opposite}, hence the name and notation;
the concept of the opposites chamber system has been introduced in
\cite{Tits:1990}).
The articles \cite{Abramenko/Muehlherr:1997}, \cite{Muehlherr}
and Lemma~6.2 of the
present paper imply the following theorem, which is a split version
of Theorem~\ref{TopCompletion} over the complex numbers.
\begin{thm}[Topological Curtis-Tits Theorem]
Consider a complex Kac-Moody group $G(A)$
associated with a symmetrizable generalized Cartan matrix
of finite size and two-spherical
type $(W,S)$.
Let $\tri^{re}$ be the set of real roots, and
$\Pi$ be a system
of fundamental roots.
Let $\I$ be the small category with objects
${\Pi\choose 1} \cup {\Pi \choose 2}$ and morphisms
$\{\alpha\} \to \{\alpha,\beta\}$,
for all $\alpha, \beta \in \Pi$, and let $\delta \colon \I \to \KOG$
be a diagram with $\delta(\{\alpha\}) = G_\alpha$,
$\delta(\{\alpha,\beta\}) = G_{\alpha\beta}$ and
$\delta(\{\alpha\} \to \{\alpha,\beta\})
= (G_\alpha \hookrightarrow G_{\alpha\beta})$.\\[2.5mm]
Then $((G(A), \tau_{KP}), (\iota_i)_{i \in { \Pi \choose 1}
\cup {\Pi \choose 2}})$, where
$\iota_i$ is the natural inclusion map, is a colimit of $\delta$
\begin{enumerate}
\item[{\rm (a)}] in the category $\G$ of abstract groups;
\item[{\rm (b)}] in the category $\HTG$ of Hausdorff topological groups;
\item[{\rm (c)}] in the category $\KOG$ of $k_\omega$ groups.
\end{enumerate}
\end{thm}
\begin{proof}
(a) follows from \cite{Abramenko/Muehlherr:1997} or \cite{Muehlherr},
(b) follows from Lemma \ref{SimpleRoots} and Part~(c) follows from~(b)
and Corollary~\ref{colimkom}.
\end{proof}
A classification of topological amalgams as for the unitary form
leading to a split analogue of Theorem~\ref{mainphan} over the complex
numbers is not immediately possible by the methods presented in this
paper, since van der Waerden's continuity theorem does not hold
for complex Lie groups (as there exist discontinuous field
automorphisms).
In view of Shtern's generalization
of this continuity theorem to arbitrary semisimple real Lie
groups~\cite{Shtern},
a classification of amalgams over the field of real
numbers is
possible once the analogue of Lemma~\ref{SimpleRoots} and
Proposition~\ref{KPfinal} has been proved for the split real form of~$G(A)$.
This leads to split analogues of Theorem~\ref{TopCompletion}
and Theorem \ref{mainphan} over the real numbers.
\section{Duality of locally {\boldmath $k_\omega$}
and almost metrizable groups}\label{secdualcat}
%
%
%
In this section, we show that the categories
of locally~$k_\omega$ abelian groups
and almost metrizable abelian groups
are dual to each other.\\[2.5mm]
Recall that a Hausdorff space~$X$
is called \emph{almost metrizable}
if each $x\in X$ is contained in a compact set~$K$
which has a countable basis
$(U_n)_{n\in \N}$
of neighbourhoods in~$X$ (see
\cite[Definition~1.22]{Aus} or \cite{Pas}).
Thus each $U_n$ is a subset of~$X$
having $K$ in its interior,
and for each subset $U\sub X$ such that
$K\sub U^0$, there exists
$n\in \N$ such that $U_n\sub U$.
It is known that a
Hausdorff topological abelian group~$G$
is almost metrizable
if and only if $G/K$ is metrizable
for a compact subgroup~$K\leq G$
(see \cite[Proposition~2.20]{Aus}~or~\cite{Pas}).\\[2.5mm]
We shall also need a related concept:
A topological space
is called \emph{\v{C}ech complete}
if it is the intersection
of a sequence of open subsets
of a compact space.
It is known that a
Hausdorff topological abelian group~$G$
is \v{C}ech complete
if and only if $G/K$ is complete and metrizable
for some compact subgroup~$K\leq G$
(see \cite[Corollary~2.21]{Aus}).
Hence every \v{C}ech complete
abelian group is almost metrizable and complete.\\[2.5mm]
The following notation will be used
in our discussions of topological
abelian groups.
Given a Hausdorff topological abelian group~$G$, we let
$\eta_G\colon G\to G^{**}$, $\eta_G(x)(\xi):=\xi(x)$
be the evaluation homomorphism
and say that $G$ is \emph{reflexive}
if $\eta_G$ is an isomorphism of
topological groups.
Given a continuous homomorphism
$f\colon G\to H$, we let
$f^*\colon H^*\to G^*$, $f^*(\xi):=\xi\circ f$
be the dual morphism.
Then $f^{**}\circ\eta_G=\eta_H\circ f$.
We write $\T_+:=\{z\in \T\colon \Repart(z)\geq 0\}$.
Given a subset $A\sub G$,
we let $A^\circ:=\{\xi\in G^*\colon \xi(A)\sub \T_+\}$
be its \emph{polar} in $G^*$.
If $B\sub G^*$,
we write ${}^\circ B:=\eta_G^{-1}(B^\circ)$
for its polar in~$G$.
Occasionally, we shall write
$\langle \xi, x\rangle:=\xi(x)$
for $x\in G$, $\xi\in G^*$.
%
%
\begin{prop}\label{dualcats}
If $G$ is an abelian locally $k_\omega$ group,
then $G^*$ is \v{C}ech complete
and hence almost metrizable and complete.
Conversely, $G^*$ is locally $k_\omega$
and complete,
for each almost metrizable abelian group~$G$.
\end{prop}
\begin{proof}
If $G$ is an abelian locally $k_\omega$ group,
then $G$ has an open subgroup $H$ which is
a $k_\omega$-group, by Proposition~\ref{insplyd}.
Then $D:=G/H$ is discrete.
Let $i\colon H\to G$ be the inclusion map and
$q\colon G\to D$ be the canonical quotient morphism.
Then $D^*$ is compact~\cite[Proposition~7.5\,(i)]{Hofmann/Morris:1998}
and $H^*$ is complete and
metrizable \cite[Propositions~2.8 and 4.11]{Aus}.
Since~$i$ is an open embedding,
$i^*$ is a quotient morphism \cite[Lemma~2.2\,(d)]{BCP}.
Furthermore, $q$ being surjective,
$q^*$ is injective and actually
an embedding since $D^*$ is compact.
By \cite[Lemma~1.4]{BCP}, we have a short exact sequence
\[
\{1\}\to
D^*\stackrel{q^*}{\to}G^*\stackrel{i^*}{\to} H^*\to\{1\}\,.
\]
Then $K:=q^*(D^*)$
is a compact subgroup of~$G^*$ and
$G^*/K\isom H^*$ is complete and metrizable.
Therefore $G^*$ is \v{C}ech complete,
using \cite[Corollary~2.21]{Aus}.\\[2.5mm]
Conversely, assume that~$G$ is almost metrizable.
Then $G/K$ is metrizable
for a compact subgroup $K\sub G$.
Let $i\colon K\to G$ be the inclusion map
and $q\colon G\to G/K=:Q$ be the canonical
quotient morphism.
Since $Q$ is metrizable,
$Q^*$ is a $k_\omega$-group \cite[Corollary~4.7]{Aus}
and complete (see \cite[Proposition~4.11]{Aus}).
By \cite[Lemma~2.5]{BCP},
$q^*\colon Q^*\to G^*$ is a continuous open homomorphism
with compact kernel. Since~$q$ is surjective,
$q^*$ is injective
and hence an open embedding.
Thus~$G^*$ has the complete $k_\omega$-group $q^*(Q^*)$
as an open subgroup. Hence~$G^*$ is complete,
and Proposition~\ref{insplyd}
shows that $G^*$ is
locally~$k_\omega$.
\end{proof}
Cf.\ \cite[Proposition~5.20]{Aus}
for a related result.
\begin{cor}\label{cechalmm}
Let $G$ be a reflexive topological abelian group.
Then
$G$ is almost metrizable
if and only if $G^*$ is locally $k_\omega$
$($in which case $G$ actually is \v{C}ech
complete, and $G^*$ is complete$)$.
Also, $G$ is locally~$k_\omega$
if and only if $G^*$ is almost metrizable
$($in which case $G$ actually is complete,
and $G^*$ is \v{C}ech complete$)$.\Punkt
\end{cor}
\begin{rem}
It is known that each abelian $k_\omega$-group
is complete~\cite{Rai}.
Hence also each abelian locally $k_\omega$ group
is complete
(exploiting Proposition~\ref{insplyd}\,(b)).
\end{rem}
\section{Dual groups of projective limits
and direct limits}\label{dirandpro}
%
%
%
The dual groups of countable direct limits of abelian
$k_\omega$-groups
and countable projective limits of metrizable abelian groups
have been studied in~\cite{ATC}.
In this section,
we describe various generalizations
of the results from \cite{ATC}.
In particular, we study the dual groups
of countable direct limits of abelian,
locally~$k_\omega$ groups and
countable projective limits of
almost metrizable abelian groups.\\[2.5mm]
Our first proposition
generalizes
\cite[Proposition~3.1]{ATC},
the proof of which given in~\cite{ATC}
requires that the limit maps~$q_i$ are
surjective.\footnote{The description of
$\pi_1(\pl\, G_\alpha)$\vspace{-.4mm}
given in the displayed formula in
\cite[proof of Proposition~3.1]{ATC}
requires that the limit map to $G_{\alpha_k}$
is surjective.}
Recall that
a subgroup $H$ of a topological abelian group
$G$ is called \emph{dually embedded}
if each character $\xi\in H^*$ extends
to a character of~$G$.
It is \emph{dually closed}
if, for each $x\in G\take H$,
there exists $\xi\in G^*$
such that $\xi|_H=1$
and $\xi(x)\not=1$.
%
%
\begin{la}\label{prop3.1}
Let $((G_i)_{i\in I}, (q_{ij})_{i\leq j})$
be a projective system of Hausdorff
topological abelian groups and continuous homomorphisms
$q_{ij}\colon G_j\to G_i$,
with projective limit $(G,(q_i)_{i\in I})$.
Then the following holds:
\begin{itemize}
\item[{\rm (a)}]
If $q_i(G)$ is dually embedded in $G_i$
for each $i\in I$ $($e.g., if each $q_i\colon G\to G_i$
has dense image$)$,
then each $\xi\in G^*$ is
of the form $\xi=\theta\circ q_i$
for some $i\in I$
and $\theta\in G_i^*$.
Hence
$G$ is dually embedded in $P:=\prod_{i\in I}G_i$.
If each $q_i$ has dense image, then furthermore
$(G^*, (q_i^*)_{i\in I})=\dl \,((G_i^*)_{i\in I},
(q_{ij}^*)_{i\leq j})$\vspace{-.8mm}
as an abstract group.
\item[{\rm (b)}]
If each $\eta_{G_i}$ is injective,
then $G$ is dually closed in $\prod_{i\in I}G_i$.
\end{itemize}
\end{la}
\begin{proof}
(a) Since $\xi^{-1}(\T_+)$
is a $0$-neighbourhood,
there exists $i\in I$ and a $0$-neighbourhood
$U\sub G_i$ such that $q_i^{-1}(U)\sub \xi^{-1}(\T_+)$
(cf.\ \cite[p.\,23]{Wei}).
Set $D:=q_i(G)$.
Since $\T^+$ does not contain any non-trivial
subgroups and $\xi(\ker(q_i))\sub \T_+$,
we deduce that $\ker(q_i)\sub \ker(\xi)$.
Hence, there exists a
homomorphism $\zeta\colon D \to \T$
such that $\zeta\circ q_i=\xi$.
Since $\zeta(U\cap D)=\zeta(q_i(q_i^{-1}(U)))
=\xi(q_i^{-1}(U))\sub\T_+$,
\cite[Lemma~2.1]{Kp2} shows
that $\zeta$ is continuous.
Now $D$ being dually embedded in~$G_i$,
the character $\zeta$ extends
to a character $\theta\colon G_i\to\T$.
Then $\theta\circ q_i=\xi$.
Furthermore, $\xi$ extends to
the continuous homomorphism
$\theta\circ \pr_i\colon P\to\T$,
where
$\pr_i\colon P\to G_i$ is the
canonical projection.
Hence $G$ is dually embedded in~$P$.
By the preceding, $G^*=\bigcup_{i\in I}q_i^*(G_i^*)$.
If we assume that each $q_i$ has dense image,
then each $q_i^*$ is injective,
entailing that $(G^*,(q_i^*)_{i\in I})$
is the asserted direct limit group.\vspace{1mm}

(b) See \cite[Lemma~5.28]{Aus}.
\end{proof}
Note that if~$I$ is countable
in the situation of Lemma~\ref{prop3.1}
and each~$q_{ij}$ is surjective,
then each $q_i$ is surjective
(as is well known).
%
%
\begin{prop}\label{thm3.2}
Let $((G_i)_{i\in I}, (q_{ij})_{i\leq j})$ be a
projective system of
reflexive topological
abelian groups,
with projective limit $(G,(q_i)_{i\in I})$.
If $q_i(G)$ is dually embedded in~$G_i$
for each $i\in I$
$($e.g., if each $q_i$ has dense image$)$,
then $\eta_G$ is open
and an isomorphism of groups.
\end{prop}
\begin{proof}
By Lemma~\ref{prop3.1},
$G$ is dually closed and dually
embedded in $\prod_{i\in I}G_i$.
Each $G_i$ being reflexive,
also $P:=\prod_{i\in I}G_i$
is reflexive (see \cite[Proposition~(14.11)]{Ban} or
\cite{Kp1}).
Now $\eta_P$ being an isomorphism of abstract groups
and open, also
$\eta_G$ is an open isomorphism of groups by
\cite[Corollary~5.25]{Aus}
(a result due to Noble).\vspace{-1mm}
\end{proof}
\begin{rem}
If $\eta_G$ happens to be continuous
in Proposition~\ref{thm3.2}
(e.g., if $G$ is a $k$-space),
then $G$ is reflexive.
Recall that subgroups of nuclear abelian groups
(as in \cite[Definition~7.1]{Ban})
are dually embedded by \cite[Corollary~8.3]{Ban}.
Hence, if a projective limit~$G$ of reflexive nuclear
groups $G_i$ has a continuous
evaluation homomorphism $\eta_G$,
then $G$ is reflexive.
See \cite{Aus}, \cite{Ban},
and \cite{Gal}
for further information
on nuclear groups.
\end{rem}
As a corollary to Proposition~\ref{thm3.2},
we obtain a generalization
of \cite[Theorem~3.2]{ATC},
where each $G_n$ was assumed a metrizable,
reflexive abelian topological group
and where each $q_{n,m}$ (and hence $q_n$)
was assumed surjective.
%
\begin{cor}\label{alsoprop}
Let $((G_n)_{n\in \N}, (q_{n,m})_{n\leq m})$ be a
projective sequence of
reflexive topological abelian groups,
with projective limit $(G,(q_n)_{n\in \N})$.
If $q_n(G)$ is dually embedded in~$G_n$
$($e.g., if $q_n$ has dense image$)$
and $G_n$ is almost metrizable
for each $n\in\N$,
then $G$ is reflexive.
\end{cor}
\begin{proof}
By Proposition~\ref{thm3.2},
$\eta_G$ is an open isomorphism of groups.
As we assume that $G_n$ is almost metrizable
and reflexive,
$G_n$ is \v{C}ech complete (by
Corollary~\ref{cechalmm}).
Thus also $G$ is \v{C}ech complete
(see \cite[Corollary~3.9.9]{Eng})
and hence a $k$-space, by 
\cite[Theorem~3.9.5]{Eng}.
Therefore $\eta_G$ is continuous
by \cite[Corollary~5.12]{Aus}
and hence an isomorphism of topological groups
(cf.\ also \cite{CaP}).
\end{proof}
Recall that a topological abelian group $G$ is called
\emph{locally quasi-convex}
if it has a\linebreak
basis of $0$-neighbourhoods $U$
such that $U={}^\circ(U^\circ)$.
Given a topological abelian group~$G$,
there is a finest locally quasi-convex
group topology $\cO_{\text{lqc}}$ on~$G$
which is coarser than the given topology;
a basis of $0$-neighbourhoods
for $\cO_{\text{lqc}}$ is given by the bipolars ${}^\circ(U^\circ)$
of $0$-neighbourhoods in~$G$
(see \cite[Proposition~6.18]{Aus}).
Then $G_{\text{lqc}}:=(G,\cO_{\text{lqc}})/\wb{\{0\}}$
is a\linebreak
locally quasi-convex group such that
each continuous homomorphism from $G$
to a Hausdorff\linebreak
locally quasi-convex group
factors over $G_{\text{lqc}}$.
In the following,
$\dl_{\text{lqc}}G_i$\vspace{-.5mm}
denotes the direct limit
of a direct system $\cS=((G_i)_{i\in I},(\lambda_{ij})_{i\geq j})$
of locally quasi-convex Hausdorff abelian groups
in the category of such groups.
It can be obtained as
\[
\dl_{\text{lqc}}G_i\; =\; (\dl\,G_i)_{\text{lqc}}\,,\vspace{-.8mm}
\]
where $\dl\,G_i=\dl\,\cS$
in the category of
Hausdorff abelian groups (cf.\ \cite[Proposition~4.2]{ATC}).\\[2.2mm]
Our next result generalizes \cite[Theorem~4.3]{ATC},
which only applied to projective sequences
$((G_n)_{n\in \N}, (q_{n,m})_{n\leq m})$
of metrizable, reflexive abelian topological
groups such that each $q_{n,m}$
(and hence $q_n$) is surjective.
%
%
%
%
\begin{prop}\label{thm4.3}
Let $\cS:=((G_i)_{i\in I},(q_{ij})_{i\leq j})$ be a projective
system of reflexive topological abelian
groups, with projective limit $(G,(q_i)_{i\in I})$.
If each~$q_i$ has dense image, then
$(G^*,(q_i^*)_{i\in I})=\dl\, ((G_i^*)_{i\in I},
(q_{ij}^*)_{i\leq j})$\vspace{-.3mm} holds
in the category of locally quasi-convex
Hausdorff topological abelian groups.
\end{prop}
\begin{proof}
Let $\cT$ be the compact-open
topology on~$G^*$.
By Lemma~\ref{prop3.1},
$G^*=\dl\, G_i^*$\vspace{-.3mm}
as an abstract group.
Since $(G^*,\cT)$ is locally quasi-convex and
the homomorphisms $q_i^*\colon G_i^*\to G^*$
are continuous, the topology $\cO$ making $G^*$
the locally quasi-convex direct limit
is finer than~$\cT$
and thus $\lambda\colon (G^*,\cO)\to (G^*,\cT)$, $\lambda(\xi):=\xi$
is continuous.
It remains to show
that every identity neighbourhood
$U$ in $H:=(G^*,\cO)$ is also an identity
neighbourhood in $(G^*,\cT)$.
We may assume that $U={}^\circ(U^\circ)$.
We write $h_i$ for $q_i^*$,
considered as a map $G_i^*\to H$;
then $h_i$ is continuous,
$\lambda\circ h_i=q_i^*$,
and $h_j\circ q_{ij}^*=h_i$ for $i\leq j$.
The continuous homomorphisms
$\eta_{G_i}^{-1}\circ h_i^*$
form a cone over the projective system~$\cS$
since $q_{ij}\circ \eta_{G_j}^{-1}\circ h_j^*
=\eta_{G_i}^{-1}\circ q_{ij}^{**}\circ h_j^*
=\eta_{G_i}^{-1}\circ h_i^*$.
Hence, there is a continuous homomorphism
$\beta\colon H^*\to G$ such that
$q_i\circ \beta=\eta_{G_i}^{-1}\circ h_i^*$
for all $i\in I$.
Then $U^\circ\sub H^*$ is compact (cf.\ \cite[Proposition~1.5]{Ban})
and hence also $K:=\beta(U^\circ)\sub G$ is compact.
We claim that $U=K^\circ$;
if we can show this, then
$U$ is also a $0$-neighbourhood in $(G^*,\cT)$
and the proof is complete.
To prove the claim, we note that
$\langle \lambda(\xi),\beta(\theta)\rangle=\langle\theta,\xi\rangle$
for each
$\theta\in H^*$ and $\xi\in H$.
In fact,
there is $i\in I$ such that
$\xi=h_i(\zeta)$ for some $\zeta\in G_i^*$.
Then
\begin{eqnarray*}
\langle\lambda(\xi),\beta(\theta)\rangle &=&
\langle \lambda(h_i(\zeta)),\beta(\theta)\rangle\;=\;
\langle q_i^*(\zeta),\beta(\theta)\rangle\;=\;
\langle \zeta, q_i(\beta(\theta))\rangle\\
&=&
\langle \zeta,\eta_{G_i}^{-1}(h_i^*(\theta))\rangle
\;=\;
\langle h_i^*(\theta),\zeta\rangle
\;=\;\langle\theta, h_i(\zeta)\rangle
\;=\; \langle\theta,\xi\rangle\,.
\end{eqnarray*}
Hence $\langle\theta,\xi\rangle=\langle \xi,\beta(\theta)\rangle$
for each $\theta\in U^\circ$ and $\xi\in G^*$
in particular,
entailing that $U={}^\circ(U^\circ)=K^\circ$.
\end{proof}
As a corollary, we obtain a generalization
of \cite[Theorem~4.4]{ATC}.
%
%
\begin{cor}\label{correfl}
Let $((G_n)_{n\in \N},(q_{n,m})_{n\leq m})$ be a projective
sequence of reflexive topological abelian
groups, with projective limit
$(G,(q_n)_{n\in \N})$.
If $q_n$ has dense image
and $G_n$ is almost metrizable
for each $n\in \N$,
then ${\dl}_{\rm lqc} G_n^*$\vspace{-.8mm}
is reflexive.
\end{cor}
\begin{proof}
By Corollary~\ref{alsoprop},
$G$ is reflexive,
whence also $G^*$ is reflexive by \cite[Proposition~5.9]{Aus}.
But $G^*=\dl_{\rm lqc}G_n^*$\vspace{-2mm}
by Proposition~\ref{thm4.3}.
\end{proof}
Part\,(a) of the next proposition generalizes
\cite[Corollary~2.2]{ATC} (the corresponding statement
for $k_\omega$-groups), while part\,(b)
generalizes \cite[Theorem~4.5]{ATC},
which applied to strict direct sequences
of abelian $k_\omega$-groups $G_n$
with $i_{n,m}(G_m)$ dually embedded in $G_n$
(cf.\ also \cite[Remark on p.\,18]{ATC}).
%
%
\begin{prop}\label{dlslcom}
Let $\cS:=((G_n)_{n\in \N}, (i_{n,m})_{n\geq m})$
be a direct sequence of abelian, locally~$k_\omega$
groups, and $(G,(i_n)_{n\in \N})$
be its direct limit in the category of Hausdorff abelian groups.
\begin{itemize}
\item[\rm (a)]
Then $(G^*, (i_n^*)_{n\in \N})
=\pl\,((G_n^*)_{n\in \N}, (i_{n,m}^*)_{n\geq m})$\vspace{-.5mm}
as a topological group.
\item[\rm (b)]
If each $G_n$
is reflexive
and each $i_n^*\colon G^*\to G_n^*$
has dense image
$($which holds, for example, if
$i_{n,m}$ is a topological
embedding and $i_{n,m}(G_m)$ is dually embedded in $G_n$
whenever $m\leq n)$,
then {\em $G_{\text{lqc}}$} is reflexive.
Furthermore,
{\em $\gamma^*\colon (G_{\text{lqc}})^* \to G^*$}
is a topological isomorphism, where
{\em $\gamma\colon G\to G_{\text{lqc}}$}
is the canonical homomorphism.
\end{itemize}
\end{prop}
\begin{proof}
(a) In the situation of (a), we have
%
\begin{equation}\label{repat}
(G^*, (i_n^*)_{n\in \N})
\;=\; \pl\,((G_n^*)_{n\in \N}, (i_{n,m}^*)_{n\geq m})\vspace{-.5mm}
\end{equation}
in the category of sets,
by the universal property
of $G=\dl\, G_n$\vspace{-.3mm}
in the category of Hausdorff groups.
Then clearly (\ref{repat})
also holds in the category of groups.
Each map $i_n^*$ being continuous,
the compact-open topology $\cO$ on~$G^*$
is finer than the projective limit topology~$\cT$.
If $K\sub G$ is compact,
then $K\sub i_n(L)$ for some $n\in \N$ and some compact subset
$L\sub G_n$ (Lemma~\ref{findco})
and thus $K^\circ\supseteq (i_n^*)^{-1}(L^\circ)$
is a $0$-neighbourhood in~$\cT$.
Thus $\cO=\cT$.

(b) By (a), we have $G^*=\pl\,G_n^*$\vspace{-.3mm}
as a topological group,
whence $G^*$ is reflexive,
by Corollary~\ref{alsoprop}.
Hence also $G^{**}$ is reflexive \cite[Proposition~5.9]{Aus}.
Since $i_n^*(G^*)$
is dense in $G_n^*$ for each $n\in\N$,
Proposition~\ref{thm4.3} shows that
$G^{**}=\dl_{\text{lqc}}\,G_n^{**}$.
Since
\[
\kappa\, :=\, \dl\,\eta_{G_n}\colon
\dl_{\text{lqc}}\,G_n\to
\dl_{\text{lqc}}\,G_n^{**}\vspace{-.6mm}
\]
is a topological isomorphism from
$G_{\text{lqc}}$ onto $G^{**}$,
we deduce that $G_{\text{lqc}}$ is reflexive.
To prove the final assertion,
note that $\kappa^*\colon G^{***}\to (G_{\text{lqc}})^*$
and $\eta_{G^*}\colon G^*\to G^{***}$
are topological isomorphisms and hence
also $\kappa^*\circ\eta_{G^*}\colon G^*\to (G_{\text{lqc}})^*$.
By the universal property of $\gamma$,
the map $\gamma^*$ is a bijection.
A simple calculation shows that
$\kappa^*\circ \eta_{G^*}=(\gamma^*)^{-1}$.
Hence $\gamma^*$ is a topological isomorphism.\\[2.5mm]
To complete the proof, note that if
$i_{n,m}$ is a topological embedding,
then each character on $G_m$
corresponds to a character
on $i_{n,m}(G_m)$,
which in turn extends
to a character of $G_n$
if we assume that $i_{n,m}(G_m)$
is dually closed in~$G_n$.
Therefore $i_{n,m}^*$
is surjective
in this case and hence also
$i_n^*$.
\end{proof}
Recall that a topological abelian group~$G$
is called \emph{strongly reflexive}
if all closed subgroups and all Hausdorff quotient
groups of~$G$ and $G^*$ are reflexive.
For example,
every \v{C}ech complete nuclear
abelian group is strongly reflexive~\cite[Theorem~20.40]{Aus}.
This enables us to strengthen
Corollary~\ref{alsoprop}
in the case of nuclear groups:
\begin{prop}
The projective limit~$G$
of a countable projective system
of \v{C}ech complete
nuclear abelian groups
is strongly reflexive.
\end{prop}
\begin{proof}
$G$ is \v{C}ech complete by
\cite[Corollary~3.9.9]{Eng}
and nuclear by
\cite[Proposition~7.7]{Ban},
whence~\cite[Theorem~20.40]{Aus}
applies.
\end{proof}
The following proposition will be
used to prove strong reflexivity of certain direct limits.
We recall that a topological abelian group~$G$
is said to be \emph{binuclear}
if $\eta_G$ is surjective and both~$G$
and~$G^*$ are nuclear~\cite[p.\,152]{Ban}.
The surjectivity of
$\eta_G$ is automatic
if $G$ is nuclear and complete
(see \cite[Theorem~6.4]{Gal}).
%
%
\begin{prop}\label{binuc}
Let $G$ be a binuclear abelian
group which is reflexive and locally~$k_\omega$.
Then $G$ is strongly reflexive.
\end{prop}
\begin{proof}
Since~$G^*$ is nuclear
and \v{C}ech complete, it
is strongly reflexive by \cite[Theorem~20.40]{Aus}.
Hence also $G\isom G^{**}$ is strongly reflexive.
\end{proof}
We first deduce a variant of
Corollary~4.6 in~\cite{ATC},
which was formulated for countable direct sums
of strongly reflexive,
nuclear abelian
$k_\omega$-groups.\footnote{The proof
of the
cited corollary given in \cite{ATC}
is problematic:
Lydia Au\ss{}enhofer~\cite{Au2} showed by example 
that the maps $(f_n^m)''$ occurring in its proof
need not be embeddings,
contrary to claims
made there.
When we wrote this article,
we did not find a way to repair
the proof if the $G_n$'s
are merely nuclear (rather than
binuclear). After the article was
completed, Lydia Au\ss{}enhofer
showed that every nuclear abelian
$k_\omega$-group
is strongly reflexive
and binuclear~\cite{LyX}.}
%
%
%
\begin{cor}\label{cor4.6}
Let $(G_n)_{n\in \N}$ be a sequence of
reflexive, binuclear abelian groups
which are locally $k_\omega$.
Then $\bigoplus_{n\in \N}G_n$
is strongly reflexive.
\end{cor}
\begin{proof}
The countable direct sum
$G:=\bigoplus_{n\in \N}G_n$
is nuclear and reflexive
by \cite[Propositions~7.8 and 14.11]{Ban},
and locally $k_\omega$ by Corollary~\ref{dsums}.
Furthermore, $G^*\isom \prod_{n\in\N}G_n^*$
by \cite[Proposition~14.11]{Ban},
whence~$G^*$ is nuclear by \cite[Proposition~7.6]{Ban}.
Thus Proposition~\ref{binuc} applies.
\end{proof}
The next corollary provides a further generalization.
%
%
\begin{cor}\label{precor4.6}
Let $((G_n)_{n\in \N}, (i_{n,m})_{n\geq m})$
be a strict direct sequence of
reflexive, binuclear abelian groups
which are locally $k_\omega$,
and $G=\dl\, G_n$\vspace{-.3mm} in the category of
Hausdorff topological abelian groups.
Then~$G$ is strongly reflexive and binuclear.
\end{cor}
\begin{proof}
Since $G$ can be realized as
a Hausdorff quotient group
of the countable direct sum
$\bigoplus_{n\in \N}G_n$,
it is nuclear by \cite[Propositions 7.5 and 7.8]{Ban}
(cf.\ also \cite[Proposition~7.9]{Ban}).
Hence~$G$ is locally quasi-convex in particular.
Furthermore, $G$ is locally $k_\omega$ (by
Corollary~\ref{beyondinj})
and reflexive, by Proposition~\ref{dlslcom}\,(b).
The latter applies because subgroups
of nuclear groups are dually embedded
\cite[Corollary~8.3]{Ban}.
By Proposition~\ref{dlslcom}\,(a),
$G^*$ is a projective limit of nuclear
groups and hence nuclear
(see \cite[Proposition~7.7]{Ban}).
Thus~$G$ is binuclear and hence
strongly reflexive, by Proposition~\ref{binuc}.
\end{proof}
\begin{rem}
It would be interesting to
find sufficient conditions
ensuring that the direct limit
Hausdorff topological group $G=\dl\,G_n$\vspace{-.5mm} of
an ascending sequence $G_1\leq G_2\leq\cdots$
of locally quasi-convex abelian $k_\omega$-groups
(or locally~$k_\omega$ groups)
is locally quasi-convex.
The only available conditions
ensuring the local quasi-convexity
of a direct limit Hausdorff topological group
$G=\dl\,G_n$\vspace{-.5mm}
seem to be the following:
\begin{itemize}
\item[(i)]
Each $G_n$ is a Hausdorff
locally convex space. Or:
\item[(ii)]
Each $G_n$ is a nuclear abelian group.
\end{itemize}
Note that
$G$ coincides with the
direct limit Hausdorff locally convex
space in the situation of~(i)
(cf.\ \cite[Proposition~3.1]{Hir}),
which is locally quasi-convex by \cite[Proposition~2.4]{Ban}).
In the situation of (ii),
$G$ is nuclear
(cf.\ \cite[Proposition~7.9]{Ban})
and hence locally quasi-convex
by \cite[Theorem~8.5]{Ban}.\\[2.5mm]
Of course, the direct limit
Hausdorff topological $G$ can always be
realized as a quotient group
of the direct sum $\bigoplus_{n\in \N}G_n$,
which is locally quasi-convex,
but quotients of locally quasi-convex
groups need not be locally quasi-convex \cite[Proposition~12.9]{Aus}.
And, as is to be expected,
there exist examples of $k_\omega$-groups
that are not locally quasi-convex~\cite{Au2},
for example the free topological vector
space over any non-discrete $k_\omega$-space
\cite[Proposition~6.4]{INS}.
\end{rem}
In the case of locally convex spaces,
no pathologies occur:
\begin{prop}\label{vec-case}
Let $E_1\sub E_2\sub\cdots$
be an ascending sequence
of locally convex real topological vector spaces~$E_n$,
such that the inclusion maps are continuous and linear.
If each $E_n$ is a $k_\omega$-space,
then the locally convex direct limit topology
on $E:=\bigcup_{n\in \N}E_n$
turns~$E$ into the direct limit
topological space $\dl\,E_n$\vspace{-.3mm},
and makes it a $k_\omega$-space.
\end{prop}
\begin{proof}
It is well known that the
box topology on $S:=\bigoplus_{n\in \N}E_n$
coincides with the locally convex direct
sum topology, and that the locally convex direct
limit $E$ can be realized as a quotient of this direct
sum by a suitable closed vector subspace.
This quotient is also the direct
limit $\dl\,E_n$\vspace{-.3mm}
in the category of Hausdorff topological groups.
Hence the assertions follow from Proposition~\ref{injsequent}.
\end{proof}
Many locally convex spaces
of interest are $k_\omega$-spaces,
e.g., all Silva spaces~\cite[Example~9.4]{DL3}.
\noindent
{\footnotesize{\bf Helge Gl\"{o}ckner} (corresponding author),
Universit\"{a}t Paderborn,
Institut f\"{u}r Mathematik,\\
Warburger Str.\ 100, 33098 Paderborn, Germany.
\,{\tt glockner@math.upb.de}\\[3mm]
{\bf Ralf Gramlich},
TU~Darmstadt, FB~Mathematik~AG~5,
Schlossgartenstr.\,7,
64289 Darmstadt, Germany.
\,{\tt gramlich@mathematik.tu-darmstadt.de}\\[2mm]
{\bf Ralf Gramlich} (alternative address),
The University of Birmingham,
School of Mathematics, Edgbaston,
Birmingham B15 2TT, United Kingdom.
\,{\tt ralfg@maths.bham.ac.uk}\\[3mm]
{\bf Tobias Hartnick}, Departement Mathematik,
HG~J~16.1,
R\"{a}mistrasse~101,
8092~Z\"{u}rich, Switzerland.
\,{\tt tobias.hartnick@math.ethz.ch}}\vfill

\begin{thebibliography}{10}\itemsep+2pt
%
%
\bibitem{Abramenko/Muehlherr:1997}
Abramenko, P.\ and B.\ M{\"u}hlherr, 
\emph{Pr{\'e}sentation des certaines $BN$-paires jumel{\'e}es comme
sommes amalgam{\'e}es},
C.\ R.\ Acd.\ Sci.\ Paris S{\'e}r.\ I Math.\ {\bf 325} (1997), 701--706.
%
%
\bibitem{AlC} Alperin, R.,
\emph{Compact groups acting
on trees}, Houston J. Math.\
{\bf 6} (1980), 439--441.
%
%
\bibitem{AlT} Alperin, R.,
\emph{Locally compact groups acting on trees},
Pacific J. Math.\ {\bf 100} (1982), 23--32.
%
%
\bibitem{ATC} Ardanza-Trevijano, S.\ and M.\,J. Chasco,
{\em The Pontryagin duality of
sequential limits of topological Abelian
groups}, J. Pure Appl.\ Algebra {\bf 202} (2005),
11--21.
%
%
\bibitem{AT2} Ardanza-Trevijano, S.\ and M.\,J. Chasco,
{\em Continuous convergence and duality of
limits of topological Abelian
groups}, in: Bol.\ Soc. Mat.\ Mexicana
{\bf 13}, (2007), no.\,1.
%
%
%
\bibitem{Aus} Au\ss{}enhofer, L.,
{\em Contributions to the duality theory of Abelian
topological groups and to the theory
of nuclear groups}, Dissertationes Math.\ {\bf 384}, 1999.
%
%
\bibitem{Au2}
Au\ss{}enhofer, L., personal communication,
March 2006.
%
%
\bibitem{LyX}
Au\ss{}enhofer, L.,
\emph{On the nuclearity of dual groups},
preprint, July 2007.
%
%
\bibitem{Ban} Banaszczyk, W., ``Additive Subgroups of Topological
Vector Spaces,'' Lecture Notes in Math.\ {\bf 1466},
Springer-Verlag, Berlin, 1991.
%
%
\bibitem{BCP}
Banaszczyk, W., M.\,J. Chasco and E. Mart\'{\i}n-Peinador,
\emph{Open subgroups and Pontryagin duality},
Math.\ Z. {\bf 215} (1994), 195--204.
%
%
\bibitem{BaB}
Beattie, R. and H.-P. Butzmann,
\emph{Continuous duality of limits and colimits
of topological Abelian groups},
Appl.\ Categor.\ Struct.\
(online first, 2007),
DOI 10.1007/s10485-007-9112-5.
%
%
\bibitem{Bennett/Gramlich/Hoffman/Shpectorov:2003}
Bennett, C.\,D., R.\ Gramlich, C.\ Hoffman and S.\ Shpectorov,
\emph{Curtis-Phan-Tits Theory}, pp.\ 13--29 in:
A.\,A.\ Ivanov, M.\,W.\ Liebeck and J.\ Saxl (eds.),
``Groups, Combinatorics and Geometry,''
2003.
%
%
\bibitem{Bennett/Shpectorov:2004}
Bennett, C.\,D. and S.\ Shpectorov,
{\em A new proof of Phan's theorem},
J.\ Group Theory {\bf 7} (2004), 287--310.
%
%
\bibitem{Benson:1991} Benson, D.\,J.,
``Representations and Cohomology. II:
Cohomology of Groups and Modules,''
Cambridge University Press, Cambridge, 1991.
%
%
\bibitem{Borovoi:1984} Borovoi, M.\,V.,
{\em Generators and relations in compact Lie groups} (Russian),
Funkts.\ Anal.\ Prilozh.\ {\bf 18} (1984), 57--58 (English translation:
Funct.\ Anal.\ Appl.\ {\bf 18} (1984), 133--134).
%
%
\bibitem{Bourbaki:2002} Bourbaki, N.,
``Lie Groups and Lie Algebras,'' Chapters 4--6,
Springer,
Berlin, 2002.
%
%
\bibitem{Brown:1989} Brown, K.\,S.,
``Buildings,'' Springer, Berlin, 1989.
%
%
\bibitem{Caprace:2005} Caprace, P.-E.,
``Abstract Homomorphisms of Split Kac-Moody Groups,''
Ph.D. Thesis, Universit\`e Libre de Bruxelles, 2005.
%
%
\bibitem{CaD} Chasco, M.\,J. and X. Dom\'{i}nguez,
\emph{Topologies on the direct sum of topological abelian
groups}, Topology Appl.\ {\bf 133} (2003), 209--223.
%
%
\bibitem{CaP} Chasco, M.\,J. and E. Mart\'{i}n Peinador,
{\em On strongly reflexive topological groups},
Appl.\ General Topol.\ {\bf 2} (2001), 219--226.
%
%
\bibitem{Curtis:1965} Curtis, C.\,W.,
{\em Central extensions of groups of Lie type},
J.\ Reine Angew.\ Math.\ {\bf 220} (1965), 174--185.
%
%
\bibitem{D} Dunlap, J.,
``Uniqueness of Curtis-Phan-Tits Amalgams,''
Ph.D.\ Thesis, Bowling Green State University, 2005.
%
%
\bibitem{Eng} Engelking, R.,
``General Topology,'' Heldermann, Berlin, 1989.
%
%
\bibitem{FaT}
Franklin, S.\,P. and B.\,V.\,S. Thomas,
\emph{A survey of $k_\omega$-spaces},
Topol.\ Proc.\ {\bf 2} (1977),
111--124.
%
%
\bibitem{Gal}
Galindo, J., {\em Structure and analysis
on nuclear groups}, Houston J. Math.\
{\bf 26} (2000), 315--334.
%
%
\bibitem{Garland:1973} Garland, H.,
{\em $p$-adic curvature and the cohomology of discrete
subgroups of $p$-adic groups}, Annals Math.\ {\bf 97} (1973), 375--423.
%
%
\bibitem{DIR} Gl\"{o}ckner, H.,
{\em Direct limit Lie groups and manifolds},
J. Math.\ Kyoto Univ.\ \textbf{43} (2003), 1--26.
%
%
\bibitem{DL2} Gl\"{o}ckner, H.,
\emph{Fundamentals of direct limit Lie theory},
Compos.\ Math.\ {\bf 141} (2005), 1551--1577.
%
%
\bibitem{DL3} Gl\"{o}ckner, H.,
\emph{Direct limits of infinite-dimensional
Lie groups compared
to direct limits in related categories},
J. Funct.\ Anal.\ {\bf 245} (2007), 19--61. 
%
%
\bibitem{INS} Gl\"{o}ckner, H.,
\emph{Instructive examples of smooth, complex differentiable
and complex analytic mappings into locally convex spaces},
J.\ Math.\ Kyoto Univ.\ {\bf 47} (2007),
631--642.
%
%
\bibitem{Goldschmidt:1980} Goldschmidt, D.,
{\em Automorphisms of trivalent graphs},
Annals Math.\ {\bf 111} (1980), 377--406.
%
%
\bibitem{Gorenstein/Lyons/Solomon:1995}
Gorenstein, D., R. Lyons and R. Solomon,
``The Classification of the Finite Simple Groups. Number~2. Part~I.
Chapter~G. General Group Theory,''
Amer.\ Math.\ Soc., Providence, 1995.
%
%
\bibitem{Gorenstein/Lyons/Solomon:1998}
Gorenstein, D., R. Lyons and R. Solomon,
``The Classification of the Finite Simple Groups. Number~3.
Part~I. Chapter~A. Almost Simple $K$-Groups,''
Amer.\ Math.\ Soc., Providence, 1998.
%
%
\bibitem{Gv2} Graev, M.\,I.,
\emph{On free products of topological groups} (Russian),
Izvestiya Akad.\ Nauk SSSR.\ Ser.\
Mat.\ {\bf 14} (1950), 343--354.
%
%
\bibitem{Grv} Graev, M.\,I.,
\emph{Free topological groups},
Amer.\ Math.\ Soc.\ Translation {\bf 1951} (1951),
no.\,35, 305--364.
%
%
\bibitem{Gramlich:2004}
Gramlich, R., ``Phan Theory,'' Habilitationsschrift, TU Darmstadt, 2004.
%
%
\bibitem{Gra} Gramlich, R.,
{\em Defining amalgams of compact Lie groups},
J. Lie Theory {\bf 16} (2006), 1--18.
%
%
\bibitem{Gramlich} Gramlich, R.,
\emph{Developments in finite Phan theory},
to appear in the Proceedings of
``Buildings and Groups'', Ghent, May 20 -- 26,
2007.
%
%
\bibitem{Gramlich/Hoffman/Shpectorov:2003}
Gramlich, R., C. Hoffman and S. Shpectorov,
\emph{A Phan-type theorem for $\Sp(2n,q)$},
J.\ Algebra {\bf 264} (2003), 358--384.
%
%
\bibitem{Gramlich/Maldeghem:2006}
Gramlich, R. and H. Van Maldeghem,
\emph{Intransitive geometries}, Proc.\ London Math.\ Soc.\ {\bf 93}
(2006), 666--692.
%
%
\bibitem{Han} Hansen, V.\,L.,
{\em Some theorems on direct limits
of expanding systems of manifolds},
Math.\ Scand.\ \textbf{29} (1971), 5--36.
%
%
\bibitem{Hartnick} Hartnick, T.,
``Moufang Topological Twin Buildings'', Diploma Thesis, Darmstadt, 2006.
%
%
\bibitem{Helgason:1978} Helgason, S.,
``Differential Geometry, Lie Groups, and Symmetric Spaces,''
Academic Press, San Diego, 1978.
%
%
\bibitem{Hilgert/Neeb:1991} Hilgert, J. and K.-H. Neeb,
``Lie-Gruppen und Lie-Algebren,'' Vieweg,
Braunschweig, 1991.
%
%
\bibitem{Hir}
Hirai, T., H. Shimomura, N. Tatsuuma and E. Hirai,
{\em Inductive limits
of topologies, their direct products, and problems
related to algebraic structures},
J. Math.\ Kyoto Univ.\ \textbf{41} (2001), 475--505.
%
%
\bibitem{Hofmann:1962} Hofmann, K.\,H.,
{\em Der Schursche Multiplikator topologischer Gruppen},
Math.\ Z.\ {\bf 79} (1962), 389--421.
%
%
\bibitem{Hofmann/Morris:1998} Hofmann, K.\,H. and S.\,A.\ Morris,
``The Structure of Compact Groups,'' De Gruyter, Berlin, 1998.
%
%
\bibitem{KacConstructingGroups} Kac, V.\,G.,
{\em Constructing groups associated to infinite-dimensional Lie
algebras} in: ``Infinite-dimensional groups with applications''
(Berkeley, Calif., 1984), Math.\ Sci.\ Res.\ Inst.\ Publ.\
\textbf{4} (1985), 167--216.
%
%
\bibitem{KacBook} Kac, V.\,G.,
``Infinite-Dimensional Lie Algebras'',
Third Edition, Cambridge University Press, Cambridge, 1990.
%
%
\bibitem{KP83} Kac, V.\,G. and D.\,H.\ Peterson 
{\em Infinite flag varieties and conjugacy theorems},
Proc.\ Nat.\ Acad.\ Sci.\ U.S.A.\ {\bf 80} (1983), 1778--1782.
%
%
\bibitem{Kac/Peterson:1985} Kac, V.\,G. and D.\,H.\ Peterson,
{\em Defining relations of certain infinite-dimensional groups},
Asterisque, 1985, 165--208.
%
%
\bibitem{Kp1} Kaplan, S.,
{\em Extensions of Pontryagin duality} I:
{\em infinite products}, Duke Math.\ J. {\bf 15} (1948),
649--658.
%
%
\bibitem{Kp2} Kaplan, S.,
{\em Extensions of Pontryagin duality} II:
{\em direct and inverse limits}, Duke Math.\ J. {\bf 17} (1950),
419--435.
%
%
\bibitem{Kat} Katz, E.,
\emph{Free products in the category
of $k_\omega$-groups},
Pacific J. Math.\ {\bf 59} (1975), 493--495.
%
%
\bibitem{KaM} Katz, E. and S.\,A. Morris,
\emph{Free products of $k_\omega$-topological groups
with normal amalgamation},
Topology Appl.\ {\bf 15} (1983), 189--196.
%
%
\bibitem{KM1} Katz, E. and S.\,A. Morris,
{\em Free products of topological
groups with amalgamation},
Pacific J. Math.\ {\bf 119} (1985),
169--180.
%
%
\bibitem{KM2} Katz, E. and S.\,A. Morris,
{\em Free products of topological
groups with amalgamation\/} II,
Pacific J. Math.\ {\bf 120} (1985),
123--130.
%
%
\bibitem{Kawada:1940} Kawada, Y.,
{\em \"Uber die \"Uberlagerungsgruppe und die stetige projektive
Darstellung topo\-logischer Gruppen},
{Jap.\ J.\ of Math.} {\bf 17} (1940), 139--164.
%
%
\bibitem{KhM} Khan, M.\,S. and S.\,A. Morris,
\emph{Free products of topological groups
with central amalgamation} I,
Trans.\ Amer.\ Math.\ Soc.\ {\bf 273} (1982), 405--416.
%
%
\bibitem{KhM2} Khan, M.\,S. and S.\,A. Morris,
\emph{Free products of topological groups
with central amalgamation} II,
Trans.\ Amer.\ Math.\ Soc.\ {\bf 273} (1982),
417--432.
%
%
\bibitem{Kramer} Kramer, L., \emph{{L}oop {G}roups and {T}win {B}uildings},
Geom.\ Dedicata {\bf 92} (2002), 145--178. 
%
\bibitem{Kumar} Kumar, S.,
``Kac-Moody Groups, their Flag Varieties and Representation Theory,''
Birkh\"auser, Basel, 2002.
%
%
\bibitem{MMO}
Mack, J., S.\,A. Morris and E.\,T. Ordman,
\emph{Free topological groups and the projective
dimension of a locally compact abelian group},
Proc.\ Amer.\ Math.\ Soc.\ {\bf 40} (1973),
303--308.
%
%
\bibitem{Mar}
Markov, A.\,A.,
\emph{On free topological groups},
Amer.\ Math.\ Soc.\ Transl.\ {\bf 1950} (1950),
no.\,30, 195--272.
%
%
\bibitem{Mitchell}
Mitchell, S.\,A., \emph{Quillen's theorem on buildings and the loops
on a symmetric space}, Enseign.\ Math.\ (2) {\bf 34} (1988), 123--166.
%
%
\bibitem{Mor}
Morris, S.\,A., \emph{Local compactness and local invariance
of free products of topological groups},
Colloq.\ Math.\ {\bf 35} (1976), 21--27.
%
%
\bibitem{MaN}
Morris, S.\,A. and P. Nickolas,
\emph{Locally compact group topologies on
an algebraic free product of groups},
J. Algebra {\bf 38} (1976), 393--397.
%
%
\bibitem{Muehlherr:1999} M{\"u}hlherr, B.,
{\em Locally split and locally finite twin buildings of $2$-spherical type},
J.\ Reine Angew.\ Math. {\bf 511} (1999), 119--143.
%
%
\bibitem{Muehlherr}
M{\"u}hlherr, B.,
{\em On the simple connectedness of a chamber system associated to a
twin building},
preprint.
%
%
\bibitem{NRW}
Natarajan, L., E.\ Rodr\'{\i}guez-Carrington
and J.\,A. Wolf,
\emph{Differentiable structure for direct limit groups},
Letters in Math.\ Phys.\ \textbf{23} (1991), 99--109.
%
%
\bibitem{NRW2}
Natarajan, L., E.\ Rodr\'{\i}guez-Carrington
and J.\,A. Wolf,
\emph{Locally convex Lie groups},
Nova J.\ Alg.\  Geom.\ \textbf{2} (1993), 59--87.
%
%
\bibitem{Num} Nummela, E.\,C.,
\emph{$K$-groups generated by $K$-spaces},
Trans.\ Amer.\ Math.\ Soc.\ {\bf 201} (1975),
279--289.
%
%
\bibitem{Ord} Ordman, E.\,T.,
\emph{Free products of
topological groups which are $k_\omega$-spaces},
Trans.\ Amer.\ Math.\ Soc.\ {\bf 191} (1974),
61--73.
%
%
\bibitem{Or3} Ordman, E.\,T.,
\emph{Free $k$-groups and free topological groups},
General Topology and\linebreak
Appl.\ {\bf 5} (1975), 205--219.
%
%
\bibitem{Pas}
Pasynkov, B.\,A., \emph{Almost-metrizable topological groups} (Russian),
Dokl.\ Akad.\ Nauk SSSR {\bf 161} (1965), 281--284
(English translation:
Soviet Math.\ Dokl.\ {\bf 7} (1966), 404--408).
%
\bibitem{PS}
Pressley, A. and Segal, G., ``Loop Groups'', Oxford Mathematical
Monographs, New York, 1986.
%
%
\bibitem{Rai}
Ra\u{\i}kov, D.\,A.,
\emph{A criterion for the completeness
of topological linear spaces and topological
abelian groups}
(Russian), Mat.\ Zametki {\bf 16} (1974), 101--106
(English translation:
Math.\ Notes.\ {\bf 16} (1974), 646--648).
%
%
\bibitem{Remy} R{\'e}my, B., {\em Groupes de {K}ac-{M}oody
d\'eploy\'es et presque d\'eploy\'es}, Ast\'erisque {\bf 277}, (2002).
%
%
\bibitem{Sch} Schubert, H.,
``Topologie,'' B.\,G. Teubner, Stuttgart, 1964.
%
%
\bibitem{Shtern} Shtern, A.\,I.,
\emph{Van der Waerden continuity theorem
for semisimple Lie groups},
Russ.\ J.\ Math.\ Phys.\ {\bf 13} (2006), 210--223.
%
%
\bibitem{Springer:1998} Springer, T.\,A.,
``Linear Algebraic Groups,'' Birkh\"auser, Basel,
1998, second edition.
%
%
\bibitem{TSH}
Tatsuuma, N., H. Shimomura, and T. Hirai,
\emph{On group topologies and unitary representations of inductive
limits of topological groups and the case of
the group of diffeomorphisms},
J. Math.\ Kyoto Univ.\ {\bf 38} (1998), 551--578.
%
%
\bibitem{Tits:1962} Tits, J.,
{\em Groupes semi-simples isotropes},
pp.\ 137--146 in: ``Colloque sur la th{\'e}orie des groupes
alg{\'e}briques,''
Bruxelles, 1962.
%
%
\bibitem{Tits:1974} Tits, J.,
``Buildings of Spherical Type and Finite $BN$-Pairs,''
Springer, Berlin, 1974.
%
%
\bibitem{Tits:1986} Tits, J.,
{\em Ensembles ordonn\'es, immeubles et sommes amalgam\'ees},
Bull.\  Soc.\  Math.\  Belg.\  S\'er.\ A {\bf 38} (1986),
367--387.
%
%
\bibitem{TitsFunctor} Tits, J.,
{\em Uniqueness and presentation of {K}ac-{M}oody groups over
fields}, J. Algebra {\bf 105} (1987), 542--573.
%
%
\bibitem{Tits:1990} Tits, J., {\em R\'esum\'e de cours}, pp.\ 87--103
in ``Annuaire du Coll\`ege de France, 89e ann\'ee, 1989--1990,'' 1990.
%
%
\bibitem{Tits:1992} Tits, J., {\em Twin buildings and groups of Kac-Moody
type}, pp.\ 249--286 in: M.W.\ Liebeck, J.\ Saxl (eds.)
``Groups, Combinatorics and Geometry,'' 1992.
%
%
\bibitem{Wei} Weil, A., ``L'integration dans les groupes
topologiques et ses applications,''
Hermann, Paris, 1965.
%
%
\end{thebibliography}
\end{document}